\documentclass[final]{siamltex}

\usepackage{epsfig}

\usepackage{amsmath,amsfonts}

\usepackage{amssymb}

\usepackage{rotating}

\usepackage{subfigure}

\newcommand{\pd}[2]{\frac{\partial#1}{\partial#2}}

\newcommand{\fu} {\mathbf{u}}

\newcommand{\fw} {\mathbf{w}}

\newcommand{\fc} {\mathbf{c}}

\newcommand{\fy} {\mathbf{y}}
\newcommand{\fz} {\mathbf{z}}
\newcommand{\fmu} {\mathbf{\mu}}

\newcommand{\fPsi} {\mathbf{\Psi}}

\title{Stochastic optimal prediction for the Kuramoto-Sivashinsky equation
        \thanks{This work was supported in part by the Applied Mathematical 
        Sciences subprogram of the Office of Energy Research of the US 
        Department of Energy under Contract DE-AC03-76-SF00098, in part 
        by the National Science Foundation under Grant DMS98-14631 and in 
        part by the Columbia University Boris A. Bakhmeteff Research 
        Fellowship in Fluid Mechanics.}}

\author{Panagiotis Stinis\thanks{Department of Applied Physics 
        and Applied Mathematics, Columbia University, New York, 
        NY 10027 ({\tt stinis@math.lbl.gov}).}}

\begin{document}

\maketitle

\begin{abstract}
We examine the problem of predicting the evolution of solutions 
of the Kuramoto-Sivashinsky equation when initial data are missing. We 
use the optimal prediction method to construct equations for the reduced
system. The resulting equations for the resolved components of the solution are 
random integrodifferential equations. The accuracy of the predictions depends
on the type of projection used in the integral term of the optimal prediction 
equations and on the choice of resolved components. The novel features of 
our work include the first application of the optimal prediction formalism to 
a nonlinear, non-Hamiltonian equation and the use of a non-invariant measure 
constructed through inference from empirical data.  
\end{abstract}



\pagestyle{myheadings}
\thispagestyle{plain}
\markboth{P. STINIS}{OPTIMAL PREDICTION  FOR THE KURAMOTO-
SIVASHINSKY EQUATION}

\section{Introduction}
Computer applications to science and engineering continue to expand.
However, there is still a wealth of problems where the computational  
power we can bring to bear is not sufficient and thus the calculations 
we perform are underresolved. This situation leads to a 
twofold question: What quantities should one strive to predict 
and what are the equations that describe the 
evolution of these quantities?

The purpose of this work is to examine the problem of prediction in the
presence of underresolution for the Kuramoto-Sivashinsky equation which models pattern 
formation in different physical contexts \cite{kur1, siv1}. We will be using the 
equation in the form
\begin{equation}
v_t+vv_x+v_{xx}+\nu v_{xxxx}=0, \label{ks}
\end{equation}
in the domain $[0,2\pi]$ with periodic boundary
conditions and initial condition $v_0(x)$. The parameter
$\nu$ has the conventional interpretation of viscosity and
controls the rate of energy dissipation at the small scales.
Depending on the value of the parameter $\nu$ the solutions can 
exhibit a wide range of behavior, from steady states,
to periodic solutions in time, to solutions with coherent spatial
structures that evolve chaotically in time (such solutions will 
be called chaotic), to solutions that exhibit temporal chaoticity
together with long range spatial statistical independence (spatio-temporal
chaos) \cite{cai,cross,mhyman2}. The problem of modelling the behavior of the 
Kuramoto-Sivashinsky equation using a reduced set of variables has proved to 
be a formidable one and different tools have been used towards this goal 
\cite{yak1,toh1,krug1,wit1}.

We are addressing the problem of underresolved computations through the method of 
optimal prediction \cite{chorin1,chorin2}. In order to proceed in the formulation of the 
equations for the quantities that will be predicted it is 
necessary to have some knowledge about the degrees of freedom that are 
unresolved. In optimal prediction what is assumed as known is an initial 
measure on the space of solutions.  What is sought is a mean solution with 
respect to this initial measure, compatible with the partial information 
available initially as well as with the limitations on the available 
computing power. This mean solution is the conditional expectation 
of the solution given the partial initial data, and is the best 
available estimate of the solution of the
full problem. The optimal prediction method consists of ways to 
construct equations that will allow the calculation of the conditional 
expectation of the solution. 

We are looking at a Galerkin 
Fourier truncation of the Kuramoto-Sivashinsky equation, thus, the solution
of the equation is determined through the integration in time of a system
of ordinary differential equations for the Fourier modes. The resolved variables 
will be a set of Fourier modes, smaller than what is 
needed for a full description of the system. The goal is to 
produce estimates for the evolution of the conditional expectations of the 
resolved Fourier modes conditioned on their prescribed initial values.

The first step towards this goal is to compute the true conditional
expectations (what will be later called the truth). To do that 
requires the knowledge of a measure on the space of all the Fourier modes 
(resolved and unresolved). The Kuramoto-Sivashinsky equation has no natural 
candidate for this measure. 
We evolve repeatedly, using initial conditions drawn from a uniform distribution, the 
equations needed for the full description of the system and collect the values of 
all the Fourier modes after some prescribed time interval.  In this way we obtain 
a collection of samples of the values of the Fourier modes. We, then, use 
maximum likelihood estimation to construct an approximation to the measure 
that describes the distribution of these samples. We construct a diagonal 
Gaussian approximation to the measure (the reason for the use of a diagonal 
measure is given later). The reason that we choose a Gaussian 
measure is that, inspection of the first few moments of the collection of samples, 
reveals that some of the modes have Gaussian statistics. Our future plan is to 
construct more elaborate approximations to the measure (e.g. mixtures of 
diagonal Gaussian measures using the Expectation-Maximization algorithm 
\cite{em}). The measure constructed is sampled, keeping the resolved 
Fourier modes fixed to their prescribed values. We use the samples as initial 
conditions to integrate the equations for all the Fourier modes. The true 
conditional expectations of  the resolved modes given their initial values, are 
the averages of the resolved modes over the samples.

The second step is to produce estimates for these conditional expectations.
To do that we need to construct equations for the evolution of the resolved 
modes. This is achieved through the use of the optimal prediction formalism
that utilizes the measure constructed through maximum likelihood estimation. 
We use the optimal prediction formalism in a certain limit called the 
short-memory approximation. This approximation is valid if the resolved modes 
depend only on their recent history. The short-memory approximation is more 
general than the delta-function approximation used frequently in the statistical 
physics literature, where the resolved modes' evolution is assumed to depend 
only on their present values. In fact, the delta-function approximation can be 
derived as a limiting case of the short-memory approximation.

The optimal prediction equations 
for the resolved modes, in the short-memory approximation, are random 
ordinary integrodifferential equations. They contain a Markovian term, 
depending only on the resolved modes' current values, a memory term that
depends on the recent values of the resolved modes and a random term
that depends on the unresolved modes.
The optimal prediction equations are solved repeatedly for 
different realizations of their random part. The initial conditions for the
resolved modes are the same for all realizations. The estimates of the 
conditional expectations of the resolved modes given their initial values 
are the averages, over the different realizations, of the values of these modes.  
We compare the estimates produced in this way to the true conditional
expectations. The results depend on the set of variables that are resolved and
on the type of projection used in the memory term. 
If the set of resolved variables includes all the linearly unstable modes and
the projection of the memory term is on the span of the resolved modes (linear 
projection), 
the agreement between the estimates of the conditional expectations
and the true conditional expectations is good for relatively long times. If the 
set of resolved variables includes all the linearly unstable modes and
the projection of the memory term is on an orthonormal set of functions of the 
resolved modes (finite-rank projection), 
the agreement between the estimates of the conditional expectations
and the true conditional expectations is good for short times only. 
If, even one linearly unstable mode is left unresolved, 
the agreement between the estimates and the true conditional expectations 
is good for short times only.

%%%%%%%%%End of Section{Introduction}%%%%%%%%%%%

\section{Optimal prediction formalism and the short-memory approximation}

We present the optimal prediction formalism in full generality and
derive from it the approximation that we will be using later, namely
the short-memory approximation. 
\subsection{Conditional expectations and the Mori-Zwanzig formalism}
Suppose we are given a system of ordinary differential equations
$$\frac{d\phi}{dt}=R(\phi)$$ with initial
condition $\phi(0)=x$, and we know only
a fraction of the initial data, say $\hat{x}$,
where $x=(\hat{x},\tilde{x})$ and correspondigly
$\phi=(\hat{\phi},\tilde{\phi})$ and that the
unresolved data are drawn from a measure
with density $f(x)$ (we will be working only with measures that are
smooth enough to have a density).

Suppose $u,v$ are functions of $x$,
and introduce the scalar product $(u,v)=E[uv]=\int{u(x)}v(x)f(x)dx$.
We will denote the space of functions $u$ with $E[u^2]< \infty$ 
by $L_2(f)$ or simply $L_2$.
We are looking for approximations of functions of $x$ by functions of
$\hat{x}$, where $\hat{x}$ are the variables that form our reduced system
(the resolved degrees of freedom). The functions of $\hat{x}$ form a
closed linear subspace of $L_2$, which we denote by $\hat{L}_2$.
Given a function $u$ in $L_2$, its conditional expectation with respect to
$\hat{x}$ is
$$E[u|\hat{x}]=\frac{\int{ufd\tilde{x}}}{\int fd\tilde{x}}.$$
The conditional expectation $E[u|\hat{x}]$ has the following properties:
\begin{itemize}

\item
$E[u|\hat{x}]$ is a function of $\hat{x}$,

\item
$E[au+bv|\hat{x}]=aE[u|\hat{x}]+bE[v|\hat{x}]$,

\item
$E[u|\hat{x}]$ is the best approximation of $u$ by a
function of $\hat{x}$:
$$E[|u-E[u|\hat{x}]|^2] \leq E[|u-h(\hat{x})|^2]$$
for all functions $h$.

\end{itemize}

We can approximate the conditional expectation $E[u|\hat{x}]$
by picking a basis in $\hat{L}_2$,
for example $h_1(\hat{x}),h_2(\hat{x}),\ldots$. For simplicity assume
that the basis functions $h_i(\hat{x})$ are orthonormal, i.e.,
$E[h_ih_j]=\delta_{ij}$. The approximation of the conditional
expectation can be written as
$E[u|\hat{x}]=\sum{a}_j h_j(\hat{x})$, where $a_j=E[u h_j]=
E[u(\hat{x},\tilde{x})h_j(\hat{x})]$. If we have a finite number
of terms only, we are projecting on a smaller subspace and the
projection is called a finite-rank projection. In the special case
where we pick as basis functions in $\hat{L}_2$ the functions
$h_1(\hat{x})=x_1,h_2(\hat{x})=x_2,\ldots,h_m(\hat{x})=x_m$,
then the corresponding finite-rank projection is called in physics
the "linear" projection (note that all projections are linear, so "linear" is used to
denote that the projection is on linear functions of the resolved variables). 
We should note here that it is not always true that $E[x_i x_j]=0$ for $i\neq{j}$.

The system of ordinary differential equations
we are asked to solve can be transformed into the linear
partial differential equation \cite{chorin4}
\begin{equation}
\label{pde}
u_t=Lu, \qquad u(x,0)=g(x)
\end{equation}
where $L=\sum_i R_i(x)\frac{\partial}{\partial{x_i}}$
is the Liouvillian and the solution of (\ref{pde}) is
given by $u(x,t)=g(\phi(x,t))$. Consider the following
initial condition for the PDE
$$g(x)=x_j \Rightarrow  u(x,t)=\phi_j(x,t)$$
We can rewrite the solution symbolically as
$$u(x,t)=e^{tL}x_j$$
This implies that for the general case where $u(x,0)=g(x)$,
the solution is $$u(x,t)=e^{tL}g(x)=g(e^{tL}x).$$ Using the
symbolic representation we can rewrite (\ref{pde}) as
$$\frac{\partial}{\partial{t}} e^{tL}x_j=L e^{tL}x_j$$
and furthermore by using the identity $L e^{tL}= e^{tL}L$
(\ref{pde}) becomes
\begin{equation}
\label{pde1}
\frac{\partial}{\partial{t}} e^{tL}x_j=e^{tL}Lx_j.
\end{equation}
If $P$ is any of the projections mentioned before and $Q=I-P$, (\ref{pde1}) 
can be rewritten as \cite{chorin4}
\begin{equation}
\label{mz}
\frac{\partial}{\partial{t}} e^{tL}x_j=
e^{tL}PLx_j+e^{tQL}QLx_j+
\int_0^t e^{(t-s)L}PLe^{sQL}QLx_jds,
\end{equation}
where we have used Dyson's formula
\begin{equation}
\label{dyson1}
e^{tL}=e^{tQL}+\int_0^t e^{(t-s)L}PLe^{sQL}ds.
\end{equation}
Equation (\ref{mz}) is the Mori-Zwanzig identity \cite{zwan1,zwan2,mori1}. 
Note that
this relation is exact and is an alternative way
of writing the original PDE. It is the starting
point of our approximations. Of course, we
have one such equation for each of the resolved
variables $\phi_j, j=1,\ldots,m$. The first term in (\ref{mz}) is
usually called Markovian since it depends only on the values of the variables
at the current instant, the second is called "noise" and the third "memory". 
The meaning of the different terms appearing in (\ref{mz}) and a connection 
(and generalization) to the fluctuation-dissipation theorems of irreversible 
statistical mechanics can be found in \cite{chorin1,stinis}.

If we write
$$e^{tQL}QLx_j=w_j,$$ 
$w_j(x,t)$ satisfies the equation
\begin{equation}
\label{ortho}
\begin{cases}
&\frac{\partial}{\partial{t}}w_j(x,t)=QLw_j(x,t) \\ 
& w_j(x,0) = QLx_j=R_j(x)-\mathfrak{R}_j(\hat{x}). 
\end{cases} 
\end{equation}
If we project (\ref{ortho})
using any of the projections discussed we get
$$P\frac{\partial}{\partial{t}}w_j(x,t)=
PQLw_j(x,t)=0,$$
since $PQ=0$. Also for the initial condition
$$Pw_j(x,0)=PQLx_j=0$$
by the same argument. Thus, the solution
of (\ref{ortho}) is at all times orthogonal
to the space of functions of $\hat{x}$. We call
(\ref{ortho}) the orthogonal dynamics equation. 

\subsection{The short-time and short-memory approximations}
The approximation we will 
examine is a short-time approximation and consists of dropping 
the integral term in Dyson's formula (\ref{dyson1}) 

\begin{equation}
e^{tQL} \cong e^{tL}  \label{sm1}.
\end{equation}
In other words we replace the flow in the orthogonal complement 
of $\hat{L}_2$ with the flow induced by the full system operator 
$L$. Some algebra shows that 
\begin{equation}
Q(e^{sQL}-e^{sL})=O(s^2) \label{sm41},
\end{equation}
and
\begin{equation}
\int_0^t e^{(t-s)L}PLe^{sQL}QLx_jds=\int_0^t e^{(t-s)L}
PLQe^{sL}QLx_jds+O(t^3). \label{sm5}
\end{equation}

As expected dropping the integral term in Dyson's formula yields an 
approximation that is good only for short times. However, 
under certain conditions this approximation can become valid 
for longer times. To see that consider the case where $P$ is 
the finite-rank projection so
\begin{equation}
PLQe^{sQL}QLx_j =\sum_{k=1}^l (LQe^{sQL}QLx_j,h_k)h_k(\hat{x}),
\label{sm6}
\end{equation}
and for the approximation
\begin{equation}
PLQe^{sL}QLx_j =\sum_{k=1}^l (LQe^{sL}QLx_j,h_k)h_k(\hat{x}).
\label{sm7}
\end{equation}

The quantities $(LQe^{sL}QLx_j,h_k)$ can be calculated from 
the full system without recourse to the orthogonal dynamics. 
Recall (\ref{sm41}) which states that the error in 
approximating $e^{sQL}$ by $e^{sL}$ is small for small $s$. 
This means that for short times we can infer the behavior of 
the quantity $(LQe^{sQL}QLx_j,h_k)$ by examining 
the behavior of  the quantity $(LQe^{sL}QLx_j,h_k)$.

If the quantities $(LQe^{sL}QLx_j,h_k)$ decay fast we can infer 
that the quantities $(LQe^{sQL}QLx_j,h_k)$ decay 
fast for short times. We cannot infer anything about the
behavior of $(LQe^{sQL}QLx_j,h_k)$ for
larger times. However, if $(LQe^{sQL}QLx_j,h_k)$ not
only decay fast initially, but, also, stay small for larger times,
then we expect our approximation to be valid for larger
times. To see this consider again the integral term in the 
Mori-Zwanzig equation. We see that the integral does 
not extend from $0$ to $t$ but only from $t-t_0$ to $t$, 
where $t_0$ is the time of decay of the quantities 
$(LQe^{sQL}QLx_j,h_k)$. 
This means that our approximation becomes

\begin{align}
\int_0^t e^{(t-s)L}PLe^{sQL}QLx_jds& \cong \int_{t-t_0}^t e^{(t-s)L}
PLQe^{sQL}QLx_jds \notag \\
&=\int_{t-t_0}^t e^{(t-s)L}PLQe^{sL}QLx_jds+\int_{t-t_0}^t 
O(s^2)ds  \notag \\
&=\int_{t-t_0}^t e^{(t-s)L}PLQe^{sL}QLx_jds+O(t^2t_0) \label{sm8}.
\end{align}
From this we conclude that the short-time approximation is 
valid for large times if $t_0$ is small and is called the 
short-memory approximation. On the other hand, if $t_0$ is 
large, then the error is $O(t^3)$ and 
the approximation is only valid for short times. Note that the
validity of the short-memory approximation can only be
checked after constructing it, since it
is based on an assumption about the large time behavior
of the unknown quantities $(LQe^{sQL}QLx_j,h_k)$. Note, that
determination of the quantities $(LQe^{sQL}QLx_j,h_k)$ requires 
the (usually very expensive) solution of the orthogonal dynamics equation.
The short-memory approximation, when valid, allows us to avoid the
solution of the orthogonal dynamics equation.

If the quantities $(LQe^{sL}QLx_j,h_k)$ do not decay
fast, then we can infer, again only for short times, that the quantities 
$(LQe^{sQL}QLx_j,h_k)$ of the exact Mori-Zwanzig equation do not
decay fast. Yet, it is possible that the quantities 
$(LQe^{sQL}QLx_j,h_k)$  start decaying very fast after short times and remain
small for longer times, 
so that the short-memory approximation could still hold. Of course, 
this can only be checked a posteriori, after the simulation of the short-memory
approximation equations.

In the statistical physics literature, the assumption that the 
correlations vanish for $s\neq 0$ is often made which is a special 
case of the short-memory approximation with the correlations replaced 
by a delta-function multiplied by the integrals. 
We shall also comment on this drastic approximation for our 
problem later when we present the numerical simulations of the short-memory 
approximation equations. An application of the short-memory approximation 
can be found in \cite{bell1}.

%%%%%%%%%End of Section{Optimal Predcition}%%%%%%%%%

\section{Density estimation}

Since the Kuramoto-Sivashinsky (KS) equation does not have a natural candidate
for the measure on the space of solutions, we have to construct an approximation
to this measure. We do so by means of maximum likelihood estimation. To 
find the maximum likelihood estimate of the parameters of the measure 
approximation we need to obtain samples of the solution of the equation. This
is accomplished through the numerical solution of the equation for different 
initial conditions sampled from a uniform distribution on the space of solutions. 
The measure constructed through maximum likelihood estimation is used later 
to compute the optimal prediction equations.

\subsection{Numerical solution of the KS equation}

As mentioned earlier, we are looking at a Fourier-Galerkin truncation of
the solution. We expand the solution in Fourier series and retain the first 
$N$ terms, 

\begin{equation*}
v_N=\sum_{k=-\frac{N}{2}}^{\frac{N}{2}-1} u_k(t) e^{ikx},
\end{equation*} 
where $u_k(t)$ are the time-dependent Fourier coefficients. 
The asymmetry between positive and negative wavenumbers comes from 
the fact that we will be working with an even number of modes. 
Since the solution of KS is real, ${u}_{-k}=u_k^{*}$, 
where, as before, $*$ denotes complex conjugation. The Fourier 
expansion transforms (\ref{ks}) into a system of ordinary differential 
equations (ODE) which reads

\begin{equation}
\frac{du_k}{dt}=-\frac{ik}{2}
\sum_{k'=-\frac{N}{2}}^{\frac{N}{2}-1}
{u}_{k'}{u}_{k-k'} + (k^2-\nu k^4)u_k \label{ksp}.
\end{equation}

The KS equation conserves the mean value $\frac{1}{2\pi}\int_0^{2\pi} 
v(x)dx=u_0$. In all subsequent calculations we assume 
$u_0=0$ without loss of generality. Linear stability analysis 
shows that the rate of growth for the mode $k$ is $k^2-\nu k^4$. 
This makes the first $[\nu^{-\frac{1}{2}}]$ (where $[]$ stands 
for the integer part) modes linearly unstable. The mode 
$[(2\nu)^{-\frac{1}{2}}]$ is the most unstable. Because of the 
presence of the fourth order stabilizing term, the stable modes 
have decay rates that increase rapidly with increasing wavenumber 
thus making the system stiff
 (see \cite{lamb1} for an extended discussion of 
 stiffness). The presence of the nonlinear coupling does 
not alter the stiffness of the ODE system. It is this property 
precisely that dictates the need for (i) an implicit solver that 
guarantees stability and (ii) a solver with variable step size to  
control the errors.

The method we chose for our numerical experiments is a 
variable-step variable-order multistep method based on 
the Backward Differentiation Formulas (BDF) that employs 
the Nordsieck representation for the storing of the necessary 
quantities at each instant \cite{hair1,gear1}. The BDF methods 
are implicit multistep methods that, up to order 6, are stable. 
Since the method of integration is implicit it results in a system of 
nonlinear algebraic equations that has to be solved iteratively and we do 
that by the Newton-Raphson method. 
Finally, the convolution sums in the RHS of the equations are evaluated 
in real space using an FFT transform and the $3/2$ rule for dealiasing 
\cite{canuto}. 
Because we are using an even number of modes the mode 
$-\frac{N}{2}$ is set equal to zero to ensure that the solution to 
KS is a real function.

\begin{figure}
\centering
\subfigure[]{\epsfig{file=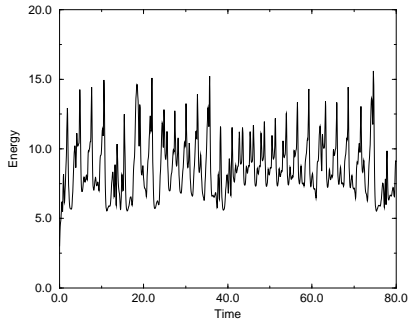,height=1.7in}}
\caption{Evolution of the energy $E=\frac{1}{4\pi}\int_0^{2\pi}
v_N^2(x,t)dx$.}
\label{fig:res1}
\end{figure}

All our calculations were performed with the value 
of the viscosity $\nu=0.085$. This value was chosen because it belongs 
in the first window that supports chaotic solutions. The first 3 Fourier modes 
are linearly unstable and mode $k=2$ is the most unstable linearly. As 
suggested by the numerical simulations, the energy spectrum shows a 
pronounced peak at this mode. Experiments with different 
number of modes show that for $N \geq 24$ the solution is converged 
for the time interval that we are interested in. For the experiments, 
an error tolerance of $10^{-7}$ per step was chosen. Simulations with 
stricter error tolerances up to $10^{-10}$ did not show any change in 
the behavior of the solutions, so we kept the value of $10^{-7}$ which 
allows for faster calculations. Strict tolerance criteria are dictated 
by the well-documented extreme sensitivity of the KS system to 
small errors\cite{mhyman1,mhyman2}. 
If high accuracy is not enforced it is possible, 
for some values of the viscosity that actually support chaotic solutions, 
to observe convergence to a steady state. On Fig.(\ref{fig:res1}) we show the 
evolution of the $L_2$ norm $E=\frac{1}{4\pi}\int_0^{2\pi}v_N^2(x,t)dx$ 
for a typical solution at this value of the viscosity. We denote the $L_2$ 
norm by $E$ and call it the energy of the system. To compute the 
maximum likelihood estimate we need to obtain a number of 
independent samples of the solution. This was done by starting the 
simulation from different random initial conditions, evolving up to time 
$t=5$ and then recording the values of the Fourier modes. The time
$t=5$ was chosen on the basis of the statistics of the collected samples.
These statistics resemble the behavior found in the literature \cite{wit1}
for long-time asymptotics. The random initial conditions were picked 
uniformly in the range $[-1,1]$ for all 
Fourier modes. We collected 10000 independent samples 
of the solution, and it is these samples that we used in the calculation of
the maximum likelihood estimate.

\subsection{Maximum likelihood estimation}

Let $\fw_1,\fw_2,\ldots,\fw_n$ denote $n$ independent 
samples of a random $d$-dimensional vector 
$\mathbf{W}$ with probability density function (p.d.f.) 
$f(\fw,\fPsi)$ where $\fPsi$ is the vector of parameters
that determine the p.d.f.. Let
\begin{equation*}
\fy=(\fw_1^T,\ldots,\fw_n^T)^T.
\end{equation*}
If we define the probability density $g(\fy,\fPsi)$ as
$$g(\fy,\fPsi)=\prod_{j=1}^n f({\fw}_j,{\fPsi}),$$
then the likelihood function formed from the observed data 
$\fy$ is defined by
\begin{align*}
L(\fPsi)=g(\fy,\fPsi).
\end{align*}
The vector $\fPsi$ is to be estimated by maximizing the likelihood \cite{lehm}. 
An estimate $\hat{\fPsi}$ of $\fPsi$ can be obtained as a solution 
of the likelihood equation 
$\pd{L(\fPsi)}{\fPsi}=0,$ or equivalently, $\pd{logL(\fPsi)}{\fPsi}=0.$

Our initial aim was to approximate the density on the space of the Fourier
modes by a mixture of Gaussian densities. Application of the maximum 
likelihood method to obtain $\fPsi,$ 
results in equations that cannot be solved in closed form. However, this
problem can be tackled by the application of iterative procedures like the
Expectation-Maximization (EM) algorithm \cite{em}. But the application of 
the EM algorithm in such a high dimensional setting as the one we have is
not straightforward. Two problems appear \cite{ormo}: The first is the problem 
of  "overfitting", where a small cluster of points, or even worse, a single 
point, determine the properties of one of the components, thus making
the entries of the covariance matrix of the component acquire 
very small values that, in their turn, create numerical  problems with the  
numerical implementation of the algorithm (we call such a component 
"narrow"). 
The second problem comes from the
fact that, if the likelihood function has multiple maxima, the EM
algorithm is only guaranteed to converge to a local maximum that
depends on the initial values used for the parameters. 
This restricts the ability of the computed mixture to be 
applied to the estimation of the density of new data points.

The KS equation is stiff and this creates an extra difficulty because
the most stable modes of the system have very small variances. 
The normalization constant of a Gaussian 
$((4\pi)^d |A|)^{-1} \exp[-\frac{1}{2}(\fw-\fmu)^{*}
A^{-1}(\fw-\fmu)$
for $d$ complex variables is $(4\pi)^{d}|A|$,
where $|A|=\prod_{i=1}^{d}(4\pi\lambda_i)$ and $\lambda_1,\ldots
\lambda_d$ are the eigenvalues of $A$. The very small variances
for the most stable modes result in very small values for the
corresponding eigenvalues, which in their turn make the normalization
constant acquire very small values. This creates, again, probelms with
the numerical implementation (see \cite{stinis}), since we have to deal with a 
"narrow" component as above. But, more importantly, it becomes increasingly 
difficult to decide, if a "narrow" component appears, whether it 
stems from overfitting or it should be there to describe the KS statistics.

Different methods have been proposed to deal with the problems
of overfitting and convergence to only local maxima for 
general Gaussian mixtures \cite{perr1,brei1,vlass1}. We shall 
adopt here a more modest goal and try to approximate the density 
with only one Gaussian. Of course, for the case
of only one Gaussian component, we can compute the properties
of the Gaussian by direct maximization of the likelihood function.
However, our future goal is to use a mixture of Gaussian components
to describe more accurately the density of the Fourier modes.

Before we estimate the parameters of the Gaussian density we make
an additional simplification dictated by the statistical properties
of the collected samples. The covariance matrix
for the Fourier modes based on the collected 
samples is
$$Cov[w_i,w_j^{*}] =\frac{1}{n-1}\sum_{k=1}^{n} (w_{ki}-M_i) 
(w_{kj}-M_j)^{*},$$ for $i,j=1,\ldots,d$ where 
$M_i, i=1,\ldots,d$ are the means
$M_i =\frac{1}{n}\sum_{k=1}^{n}w_{ki}.$ 
Our numerical experiments indicate that the covariance 
matrix for the KS system is almost diagonal, 
i.e. the non-diagonal entries are very small compared to 
the diagonal ones (results from \cite{kevin} also support this 
conclusion). We discard the off-diagonal elements and work with a Gaussian
density with a diagonal matrix.

The fact that the solution of the KS equation is real induces a constraint 
on the values of the Fourier modes, namely
for each wavenumber $k$ we have, $u_{-k}=u_k^{*}$. We can incorporate 
the constraint for the Fourier modes in the
expression for the Gaussian density. If we order the Fourier modes
so that the first $\frac{d}{2}$ positions of the vector $\fw$ of modes 
are occupied by the modes with positive wavenumbers 
$1,\ldots,\frac{d}{2}$ and the rest $\frac{d}{2}$ positions 
by the modes with negative wavenumbers $-1,\ldots,\,-\frac{d}{2},$
we can write the Gaussian density (with the constraint) as
\begin{equation}
Z^{-1} \exp[-\frac{1}{2}(\fw-\fmu)^{*}
A^{-1}(\fw-\fmu)] \prod_{i=1}^{\frac{d}{2}} 
\delta(w_{i+\frac{d}{2}}-w_i^{*}),
\label{gm13}
\end{equation}
where $A=diag(a_{1},\ldots,a_{d}).$ The entries of the matrix
$A$ have the property $a_{i+\frac{d}{2}}=a_i$ for 
$i=1,\ldots,\frac{d}{2},$ since they are the variances of the
modes $u_i$ and $u_{i+\frac{d}{2}}$ which according to our
ordering are complex conjugates of each other. Also, 
$\mu_{i+\frac{d}{2}}=\mu_i^{*}$ for $i=1,\ldots,\frac{d}{2}.$ 
The normalization constant $Z$ (modified due to the constraint) is 
$Z=\prod_{i=1}^{\frac{d}{2}}(2\pi a_i).$

The numerical results of the simulations for the KS equation are 
converged with $N=24$ modes. However, because the zero mode 
(i.e. the mean value of the real-space solution) is constant and taken equal 
to zero and the $-\frac{N}{2}$th mode was set equal to zero as 
explained above, the effective number of modes is $N=22$. Thus 
the dimension $d$ for the Gaussian density is 22. The formulas for the 
estimation of the parameters of the Gaussian density are 
\begin{align}
\mu_i&= \frac{1}{n} \sum_{j=1}^{n}  w_{ji}
 \label{gm9}, \\ 
A_{ii} &=\frac{1}{n} \sum_{j=1}^{n}
(w_{ji}-\mu_i)(w_{ji}-\mu_i)^{*} 
\label{gm10},
\end{align}
for $i=1,\ldots,d.$ In order to check how well the 
Gaussian describes the statistics of the collection of samples, we 
computed the first 4 moments (mean, covariance, skewness and flatness) 
of the samples. The Gaussian density determined by maximum likelihood 
reproduces very accurately the first three moments of the 
samples. The third moment of a Gaussian density is 
identically zero and thus it means that the distribution 
of a Fourier mode is symmetric around its mean value. 
The KS dynamics become evident in the fourth moment. 
The 4th order experimental statistics show a strong deviation from 
Gaussianity that manifests itself as negative flatness for the most 
unstable modes and positive flatness for the most stable modes. 
A negative flatness for the most unstable modes suggests the existence
of  large weights for values close to the mean value. Such large weights
in the p.d.f. of the most unstable modes have been related 
(\cite{wit1}) to the appearance of coherent structures at these scales. 
On the other hand, a positive flatness for the most stable 
modes suggests the existence of large weights at the tails of the p.d.f.. 
They have been related to the appearance of rare bursts of activity
at these scales. These rare bursts of activity are the space-time
defects (creation and annihilation of basic 
cell-structures) mentioned earlier.

The fact that a single Gaussian density cannot model accurately 
all the moments means that this density is not invariant for the 
equation. This will result in a complication for the formulas of 
the optimal prediction equations. If the density was 
invariant and we used the finite-rank projection we could have a 
fluctuation-dissipation relation \cite{stinis} which would yield significantly 
simpler expressions for certain terms of the optimal prediction 
equations.

\begin{figure}
\centering
\subfigure[]{\epsfig{file=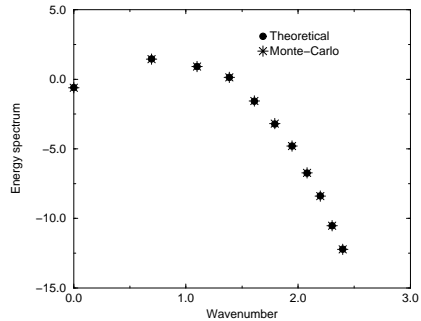,height=1.7in}}
\caption{Log-log plot of the variances for different Fourier modes as
computed by 10000 Monte-Carlo samples and as predicted by the 
Gaussian density.}
\label{fig:res2}
\end{figure}

The Gaussian density constructed through maximum likelihood can be
sampled via Monte-Carlo \cite{binder}, or by producing independent samples 
of Gaussian distributions of the appropriate variance for each of the Fourier 
modes (e.g. using the Box-Mueller algorithm) \cite{kloeden}. The results 
can be seen in Fig.(\ref{fig:res2}).

%%%%%%%%%End of Section{Density estimation}%%%%%%%%%%%

\section{Optimal prediction for Kuramoto-Sivashinsky}

We proceed to construct the short-memory approximation optimal prediction 
equations for KS. To conform with the optimal prediction formalism we set
\begin{equation*}
R_k({\fu})=-\frac{ik}{2}\sum_{k'=-\frac{N}{2}}^{\frac{N}{2}-1}
u_{k'}u_{k-k'} + (k^2-\nu k^4)u_{k}
\end{equation*}
and we get
\begin{equation}
\frac{du_k}{dt}=R_k({\fu})
\label{kspsm0}
\end{equation}
for $k=-\frac{N}{2},\ldots,\frac{N}{2}-1$. The system
of equations is supplemented with the initial condition 
$\fu(0)=\fu_0$. The initial condition can be written as
$\fu_0=(\hat{\fu}_0,\tilde{\fu}_0)$, where 
$\hat{\fu}_0$ is the vector of resolved variables
and $\tilde{\fu}_0$ the vector of unresolved ones.
Recall, that the KS equation conserves the zero
mode $u_0$ which we have set to zero for
our numerical experiments. Also, 
the $-\frac{N}{2}$th mode is set to zero to preserve
the reality of the solution of KS. This makes the effective
number of modes for the truncation equal to $N-2$. This
$(N-2)$-dimensional vector is the vector of all the modes,
resolved and unresolved. We set $n=N-2$ to conform with the 
optimal prediction formalism. In what follows
we will refer to the system of $n$ Fourier modes as the full system.
The vector of resolved variables will be
a susbset of the $n$-dimensional vector of Fourier modes 
containing $m$ modes. The unresolved modes are the rest
$n-m$ modes. The $m$-dimensional vector of resolved modes
contains $\frac{m}{2}$ positive wavenumber modes and their
corresponding $\frac{m}{2}$ negative wavenumber modes.
Accordingly, the vector of unresolved modes contains
$\frac{n-m}{2}$ positive wavenumber modes and their
corresponding $\frac{n-m}{2}$ negative wavenumber modes.
We do not specify which Fourier modes will be resolved 
since in our numerical experiments we use different sets 
of Fourier modes as resolved. By relabeling the resolved Fourier
modes we can order them from $1$ to $m$ and thus
we can construct the optimal prediction equations without
any explicit reference to the set of resolved modes.

The equations for the evolution of the set of
$m$ resolved variables in the short-memory approximation are
\begin{equation}
\frac{\partial}{\partial{t}} e^{tL}u_{0j}=
e^{tL}PLu_{0j}+e^{tL}QLu_{0j}+
\int_0^t e^{(t-s)L}PLQe^{sL}QLu_{0j}ds,
\label{kspsm}
\end{equation}
where $j=1,\ldots,m$, $L=\sum_{i=1}^{n}R_i({\fu}_0)
\pd{}{u_{0i}}$ and $u_{0j}=u_j(0)$. Since 
$PLu_{0j}=\mathfrak{R}_j(\hat{\fu}_0)$, $Lu_{0j}=R_j({{\fu}}_0)$ 
we have
\begin{equation}
QLu_{0j}=R_j({{\fu}}_0)-\mathfrak{R}_j(\hat{\fu}_0).
\label{kspsm00}
\end{equation}

\subsection{Markovian term}

We proceed with the calculation of the first term  
$PLu_{0j}=\mathfrak{R}_j(\hat{\fu}_0)$ in (\ref{kspsm}). For the case of a 
Gaussian density the conditional expectations can be computed explicitly 
\cite{chorin3}.

Let $G$ be the 
$m\times n$ matrix that has the property $G\fu_0=\hat{\fu}_0$. For 
the diagonal Gaussian density we find 
\begin{equation}
E[u_{0i}|\hat{\fu}_0]=\sum_{k=1}^{m} q_{ik} u_{0k} +c_i \label{ce}
\end{equation}
for $i=1,\ldots,n$ where $q_{ik}$ are the entries of the $n\times m$
matrix $Q$ and $c_i$ are the entries of the $n$-vector $\fc$ given
respectively by 
\begin{align}
Q &= (CG^{*})(GCG^{*})^{-1} \label{ceq} \\
\fc &=\fmu-(CG^{*})(GCG^{*})^{-1}(G\fmu) \label{cec},
\end{align}
where $C$ is the covariance matrix of the Gaussian density 
with entries $C_{jl}=E[(u_{0j}-\mu_j)(u_{0l}-\mu_l)^{*}]=A_{jl}$
and $\fmu$ is the vector of expectation values. Note that in the
case of a diagonal Gaussian density, the entries 
$q_{ik}$ for $i=m+1,\ldots,n$ and $k=1,\ldots,m$ are zero. Hence,
the conditional expectations of the unresolved modes conditioned
on the resolved ones are equal to their unconditional expectations. 
This is as expected, because if the Gaussian is diagonal, the modes are
independent from each other.

The conditional covariance matrix has entries 
\begin{equation}
Cov[u_{0i},u_{0j}^{*}]_{\hat{\fu}_0} =
E[u_{0i} u_{0j}^{*}|\hat{\fu}_0]-E[u_{0i}|\hat{\fu}_0]
E[u_{0j}^{*}|\hat{\fu}_0]
 = [C-QGC]_{ij} \label{cev},
 \end{equation}
for $i,j=1,\ldots,n$. For the case of a diagonal Gaussian density,
the entries of the conditional covariance matrix corresponding
to interactions among resolved modes and their interaction with the 
unresolved modes are zero, for the same reasons explained above. 
Also, the entries of the covariance matrix for the interactions among 
unresolved modes are equal to the entries of the 
unconditional covariance matrix.

Wick's theorem holds for conditional 
expectations too \cite{chorin3}. 
Since we are dealing with complex variables, 
higher moments than the second consist of products of 
centered variables $u_{0i}-E[u_{0i}|\hat{\fu}_0]$ that have to be 
defined in a way that for even moments 
the resulting product is a positive number. For this reason 
we define the variable $w_i=u_{0i}$ or 
$w_i=u_{0i}^{*}$ depending on the product and find 
\begin{equation}
E[\prod_{p=1}^{P} r_p|\hat{\fu}_0]=
\begin{cases}
0 & \text{P odd}, \\
\sum_{perm} 
Cov[w_{i_1},w_{i_2}^{*}]_{\hat{\fu}_0}\ldots 
Cov[w_{i_1},w_{i_2}^{*}]_{\hat{\fu}_0}
& \text{P even},
\end{cases}
\end{equation}
where the summation is over all the possible pairings of the P indices
and $r_p=w_p-E[w_p|\hat{\fu}_0]$.

The nonlinear sum in (\ref{kspsm0}) can be decomposed into 3 parts as
\begin{equation}
\sum_{<_1,<_2}u_{<_1}u_{<_2}
+\sum_{<_1,>_2}u_{<_1}u_{>_2}
+\sum_{>_1,>_2}u_{>_1}u_{>_2} \label{kspd},
\end{equation}
where $\sum_{<_1,<_2}$ indicates that we sum over 
all pairs of indices where both indices belong in the resolved 
range $i=1,\ldots,m$, $\sum_{<_1,>_2}$ indicates summation 
over all pairs of indices where only one index is in the resolved 
range and $\sum_{>_1,>_2}$ indicates summation over the pairs where 
both indices are outside the resolved range. 
 
Using the decomposition above can facilitate the computation of the
conditional expectations. We have
\begin{align}
E[R_j({\fu_0})|\hat{\fu}_0] &=-\frac{ij}{2}
\sum_{<_1,<_2}u_{0<_1}u_{0<_2} 
-\frac{ij}{2}\sum_{<_1,>_2}u_{0<_1}E[u_{0>_2}|\hat{\fu}_0] 
\label{kspfop1} \\
& \quad -\frac{ij}{2}\sum_{>_1,>_2}
E[u_{0>_1}u_{0>_2}|\hat{\fu}_0]  
+(j^2-\nu j^4)u_{0j} \notag
\end{align}
For the quantity $E[u_{0>_2}|\hat{\fu}_0]$ we have
\begin{equation}
E[u_{0>_2}|\hat{\fu}_0]=\sum_{l=1}^{m}q_{>_2,l}u_{0l}+c_{>_2},
\label{kspfop2}
\end{equation}
while for the quantity $E[u_{0>_1}u_{0>_2}|\hat{\fu}_0]$ we have
\begin{align}
E[u_{0>_1} u_{0>_2}|\hat{\fu}_0] &=
E[u_{0>_1} u_{0->_2}^{*}|\hat{\fu}_0] \notag \\
&= [C-QGC]_{0>_1,0->_2}+E[u_{0>_1}|\hat{\fu}_0]
E[u_{0->2}^{*}|\hat{\fu}_0].
\label{kspfop3}
\end{align}
The expressions (\ref{kspfop1}),(\ref{kspfop2}) and (\ref{kspfop3}) 
determine the Markovian term. If we use a diagonal Gaussian density, 
some of the terms in these expressions become zero. However, we kept
the most general form for reasons of completeness of the presentation.

\subsection{Noise term}

The second term in (\ref{kspsm}) is the noise term
$QLu_{0j}$, and we proceed with its calculation.
 From (\ref{kspsm00}), (\ref{kspd}) and 
(\ref{kspfop1}) we have
\begin{align}
QLu_{0j} &= -\frac{ij}{2}\sum_{<_1,>_2}u_{0<_1}
(u_{0>_2}-E[u_{0>_2}|\hat{\fu}_0] ) \label{kspsm1} \\
& \quad -\frac{ij}{2}\sum_{>_1,>_2}
(u_{0>_1}u_{0>_2}-
E[u_{0>_1}u_{0>_2}|\hat{\fu}_0]) \notag \\
&=A_j({\fu}_0)+B_j({\fu}_0), \notag  
\end{align}
where  $$A_j({\fu}_0)=-\frac{ij}{2}
\sum_{<_1,>_2}u_{0<_1}
(u_{0>_2}-E[u_{0>_2}|\hat{\fu}_0] )$$ and
$$B_j({\fu}_0)=-\frac{ij}{2}\sum_{>_1,>_2} 
(u_{0>_1}u_{0>_2}- 
E[u_{0>_1}u_{0>_2}|\hat{\fu}_0]).$$
Defining similarly $A_j({\fu}({\fu}_0,s))$ and 
$B_j({\fu}({\fu}_0,s))$ such that
$$e^{sL}QLu_{0j}=A_j({\fu}({\fu}_0,s))+
B_j({\fu}({\fu}_0,s))$$ 
we have in the short-memory approximation 
\begin{align}
\frac{\partial}{\partial{t}} e^{tL}u_{0j} &=
e^{tL}\mathfrak{R}_j(\hat{\fu}_0) \label{kspsm2} \\
& \quad +A_j({\fu}({\fu}_0,t))+
\int_0^t e^{(t-s)L}PLQA_j({\fu}({\fu}_0,s))ds \notag \\
& \quad +B_j({\fu}({\fu}_0,t))+
\int_0^t e^{(t-s)L}PLQB_j({\fu}({\fu}_0,s))ds \notag
\end{align}
From (\ref{kspsm1}) and (\ref{kspsm2}) we see
that each one of the sums $\sum_{<_1,>_2}$ 
and $\sum_{>_1,>_2}$ make a contribution 
to the memory and the noise. We should mention here that the
emergence and significance of multiplicative noise terms like 
the one coming from the sum $\sum_{<_1,>_2}$ have been also discussed in
the context of the stochastic modelling reduction strategy of Majda {\it et al.} 
\cite{majda1,majda2}.

For the choices that we 
will make for the set of resolved variables, we can safely 
remove the contribution to the memory and the noise that 
comes from $\sum_{>_1,>_2}$ as small due to the stiffness
of the system; the stiffness of the KS system is due to the
large ratio of growth rates between the small wavenumber modes
and the large wavenumber modes that are present in the solution. 
In all the numerical experiments, the set of resolved variables includes 
mostly small wavenumber modes, thus the bulk of the unresolved modes 
are large wavenumber modes. These large wavenumber modes have 
small magnitude and thus the contribution of their interaction 
to the memory and noise through the term $\sum_{>_1,>_2}$
is also small and thus can be neglected. However, recall that we 
retain the contribution of  $\sum_{>_1,>_2}$  to the Markovian (first) term
$\mathfrak{R}_j(\hat{\fu}_0)$.

The expression for $A_j({\fu}({\fu}_0,t))$ involves quantities that are 
unknown, namely the modes $u_{>_2}(t)$. We will approximate these modes 
as stochastic processes with known mean and autocorrelation by using 
the moving average method \cite{gikh1,mcco1,yagl1,chris}. Due 
to the fact that we are using a density that is non-invariant, the 
evolution of all the Fourier modes, and of the unresolved
in particular, does not constitute a stationary stochastic process. 
However, for our numerical experiments we will treat the unresolved
modes' evolution as a stationary stochastic process and use their statistics 
with respect to the initial density for the moving average method simulations. 
We defer the application of more sophisticated sample path generating methods 
suitable for nonstationary processes for future work.

Let $u(t,\omega)$ be a wide sense stationary stochastic process. If its covariance 
$$R(t_1,t_2)=E[(u(t_2,\omega)-m(t_2))(u(t_1,\omega)-m(t_1))^{*}]$$ is 
written as $R(t_1,t_2)=\int e^{ik(t_2-t_1)}dF(k)$ for some non-decreasing 
function $F(k)$ (through Khinchin's theorem) and also $F(k)$ is such that
$dF(k)=\phi(k)dk$, then 
\begin{equation}
u(t,\omega) = \int h^{*}(s-t)\rho(ds)
\label{mam}
\end{equation}
where the function $h(t)$ is the inverse Fourier transform of 
$\hat{h}(k)=\sqrt{\phi(k)}$ and also, $E[|\rho(ds)|^2]=ds.$ Note that the 
random measure $\rho$ constructed as increments of Brownian motion at instants 
$ds$ apart has this property. 
Thus, any wide-sense stationary stochastic process with $dF(k) =\phi(k)dk$
can be approximated as a sum of translates (in time) of a fixed function,
each translate multiplied by independent Gaussian random variables.
This is the "moving average" representation.

Having developed the formalism for the moving average representation,
we now derive formulas for its numerical implementation. 
A possible approach would be to discretize the expression (\ref{mam})
in intervals of length $\Delta t$ and use a midpoint rule to find
\begin{equation}
u(t_j)=\sum_{i=-n+j}^{n+j} h^{*}(t_i-t_j)\rho_i
\label{mamd}
\end{equation}
with $E[|\rho_i|^2]=\Delta t$ and $t_j=j \Delta t$. 
(\ref{mamd}) can be rewritten (by the change of variables 
$i'=i-j$ and dropping the primes) as
\begin{equation}
u(t_j)=\sum_{i=-n}^{n} h^{*}(t_i)\rho_{i+j}
\label{mamd2}
\end{equation}
and in effect the velocity at $t_j$ is the dot product of
a vector of values of $h$ and a random vector. Recalling the remark
above about the increments of Brownian motion, we see that
the process $u$ can be represented as a sum of suitably weighted
Gaussian independent random variables with mean zero and variance
$\Delta t$. In the numerical experiments, we split each unresolved mode
in its real and imaginary parts which are approximated independently.

\subsection{Memory term}

After computing the Markovian and noise terms, we proceed
with the calculation of the memory term 
\begin{equation}
\int_0^t e^{(t-s)L}PLQA_j({\fu}({\fu}_0,s))ds. 
\label{kspsm3}
\end{equation}
Note that as in the case of the noise we only use the
part of the memory that comes from the $\sum_{<_1,>_2}$ 
terms in the convolution sums.
We will use two different projections $P$, namely the
linear and the finite-rank one. The conditional expectation
projection will not be used because lack of explicit 
expressions make it computationally very expensive.

\subsubsection{Linear projection}

For a function $g({\fu}_0)$ the linear projection is
\begin{equation}
(Pg)(\hat{\fu}_0)=\sum_{i,j=1}^{m} b_{ij}^{-1}(g,u_{0i})u_{0j}
\label{linear}
\end{equation}
where $b_{ij}^{-1}$ are the entries of a matrix whose inverse
has entries $b_{ij}=(u_{0i},u_{0j})^{*}$ and the inner product
is defined through the density $f(x)$  (Gaussian in
our case) as
$$(g,h)=\int_{-\infty}^{+\infty}g(x)h^{*}(x)f(x)dxdx^{*}.$$
Using (\ref{linear}) we have
\begin{align}
\int_0^t e^{(t-s)L}PLQA_j&({\fu}({\fu}_0,s))ds = 
\label{kspsm4}\\
& \quad \int_0^t e^{(t-s)L} \sum_{i,k=1}^{m}b_{ik}^{-1}
(LQA_j({\fu}({\fu}_0,s)),u_{0i})u_{0k}ds. \notag
\end{align}
What remains to be done is to compute
the quantity $LQA_j({\fu}({\fu}_0,s))$. Some algebra and use 
of the relation
$$\sum_{l=1}^{n}R_l({\fu}_0) \pd{}{u_{0l}}
A_j({\fu}({\fu}_0,s))=
\sum_{l=1}^{n}R_l({\fu}({\fu}_0,s))
(\pd{A_j}{u_{0l}})
({\fu}({\fu}_0,s)),$$
gives
\begin{align}
L&QA_j({\fu}({\fu}_0,s))= \label{kspsm5} \\
&\sum_{l=1}^{n}R_l({\fu}({\fu}_0,s))
(\pd{A_j}{u_{0l}})
({\fu}({\fu}_0,s)) -
\sum_{i,k=1}^{m}b_{ik}^{-1}
(A_j({\fu}({\fu}_0,s)),u_{0i})R_k({\fu}_0). \notag
\end{align}
Then the 
quantities $(LQA_j({\fu}({\fu}_0,s)),u_{0i})$ can be computed 
by sampling the density via Monte Carlo, 
evolving the full system and averaging.

\subsubsection{Finite-rank projection}

The finite-rank projection of a function $g({\fu}_0)$ 
on a finite number of terms of  an orthonormal set of functions 
$h_1(\hat{\fu}_0),h_2(\hat{\fu}_0),\ldots$ is
\begin{equation}
(Pg)(\hat{\fu}_0)=\sum_{i=1}^{l} (g,h_i)h_i,
\label{finite}
\end{equation}
where the inner product is defined as before and 
$(h_i,h_j)=\delta_{ij}$. Due to the fact that we are using 
a complex Gaussian density we can define an orthonormal set of
functions
by suitably modifying one-dimensional Hermite polynomials  
for each of the $m$ resolved variables. Let 
$\hat{\fu}_0=(u_{01},\ldots,u_{0m})$ be the vector of initial values
of the resolved variables, where $u_{0j}=z_{j1}+i z_{j2}$.
Define the $2m$-dimensional  real vectors 
$$\fz=(z_{11},z_{12},\ldots,z_{m1},z_{m2})$$
and
$$\fmu=(\mu_{11},\mu_{12},\ldots,\mu_{m1},\mu_{m2}).$$
We need to split each mode in its real and imaginary parts
because of the definition of the orthonormal set below.

We define the multi-index $\kappa=(\kappa_1,\ldots,\kappa_m)$ and 
the set of functions
\begin{equation}
h^{\kappa}(\fz)=
\prod_{j=1}^{m}[\tilde{H}_{\kappa_j}
(\sqrt{\frac{2}{a_j}}
(z_{j1}-\mu_{j1}))+i \tilde{H}_{\kappa_j}(\sqrt{\frac{2}{a_j}}
(z_{j2}-\mu_{j2}))],
\label{finite1}
\end{equation}
where 
\begin{equation}
\tilde{H}_{\kappa_j}(x)=\frac{1}{\sqrt{2}}
(1+2\beta_{\kappa_j})^{\frac{1}{4}}
H_{\kappa_j}((1+2\beta_{\kappa_j})^{\frac{1}{2}}x)
e^{-\frac{\beta_{\kappa_j}}{2}x^2}.
\label{finite2}
\end{equation}
The functions $H_{\kappa_j}$ are Hermite polynomials (with weight
$\exp(-\frac{1}{2}x^2)$) satisfying
\begin{equation}
H_0(x)=1, \quad H_1(x)=x, \quad H_k(x)=\frac{1}{\sqrt{k}}
xH_{k-1}(x)-\sqrt{\frac{k-1}{k}}H_{k-2}(x)
\label{finite3}
\end{equation}
and $\beta_{\kappa_j}> -\frac{1}{2}$. The derivatives of the functions 
$\tilde{H}_{\kappa_j}(x)$ can be computed by the recursive relation
\begin{equation}
\frac{d}{dx}\tilde{H}_{\kappa_j}(x)=(1+2\beta_{\kappa_j})^{\frac{1}{2}}
\tilde{H}_{\kappa_{j-1}}(x)-\beta_{\kappa_j} x \tilde{H}_{\kappa_j}(x).
\label{finite4}
\end{equation}
To ensure the orthonormality of the above set of functions, we have to impose a 
constraint on the values of the entries of the multi-index $\kappa.$ The 
constraint is, that the entries corresponding to a positive wavenumber and its 
corresponding negative wavenumber cannot be simultaneously different from 
zero. We, also, have to enforce the constraint $\beta_{\kappa_j}=0$ when 
$H_{\kappa_j}$ is of order zero. These constraints arise from the fact that
the solution of Kuramoto-Sivashinsky equation is real thus the $k$-th negative
Fourier mode is the complex conjugate of the $k$-th positive mode.

The memory term becomes
\begin{align*}
\int_0^t e^{(t-s)L}PLQA_j&({\fu}({\fu}_0,s))ds = \\
& \quad \int_0^t e^{(t-s)L} \sum_{\kappa \in I}
(LQA_j({\fu}({\fu}_0,s)),h^{\kappa}(\hat{\fu}_0))
h^{\kappa}(\hat{\fu}_0)ds,
\end{align*}
where
\begin{align*}
LQA_j({\fu}&({\fu}_0,s))=\\
&\sum_{i=1}^{n}R_i({\fu}({\fu}_0,s))
(\pd{A_j}{u_{0i}})
({\fu}({\fu}_0,s)) \\
& - \sum_{i=1}^{m}R_i({\fu}_0)
\sum_{\kappa \in I}
(A_j({\fu}({\fu}_0,s)),h^{\kappa}(\hat{\fu}_0))
\pd{h^{\kappa}(\hat{\fu}_0)}{u_{0i}}
\end{align*}
and $I$ denotes the set of $m$-tuples of indices
used in the finite-rank projection.

Thus, we see that starting from the 
short-memory approximation we have transformed the original system 
of ODEs into a system of random ordinary integrodifferential equations. 
Such an approach was adopted in \cite{bell1} and termed stochastic 
optimal prediction. In \cite{bell1}, the form of the projection coefficients 
allowed the reduction of the integrodifferential equations to differential
equations (this is the delta-function approximation dicussed earlier).

%%%%%%%%%End of Section{Optimal prediction for KS}%%%%%%%%

\section{Numerical simulations}

We use the short-memory approximation equations
for the KS system to compute the evolution of the conditional expectations 
of a set of resolved Fourier modes conditioned on their initial values.  
We conducted numerical experiments for two sets of resolved variables.
The first set includes all the modes that are linearly unstable (for the value 
of the viscosity used). The second set includes all but one linearly unstable
modes. The error exhibited by the short-memory approximation is 
strikingly different for these two sets of resolved variables and we offer an 
explanation for this difference.

\subsection{Resolution of all unstable modes}

We present the results of the short-memory approximation for the
first set of resolved variables that includes all the unstable modes.

\subsubsection{Linear projection}

In order to construct the short-memory approximation equations
we have to compute the first (Markovian), second (noise) and 
third (memory) terms in the RHS of the short-memory approximation
equation. The Markovian term can be computed
explicitly, while the noise and memory terms rely on expressions
that have to be computed through simulations of the full system.

For the value of viscosity 
$\nu=0.085$ used in our experiments, a truncation retaining the first 
$N=24$ modes is enough to fully resolve the system. Since two of
the modes are set to zero, the effective number of modes is $n=N-2=22.$ 
The $n$-dimensional system is called the full system. We assume that initially 
we know the values of only $m$ of the $n$ modes, while the values of the 
rest $n-m$ modes are drawn from the diagonal Gaussian density (\ref{gm13}). 
The $m$ resolved modes constitute the 
reduced system. We examine the case where we know initially 
only $\frac{N}{2}=12$ variables. However, for the same reasons as in the 
case of the full system,  the zeroth and -6th mode are set to zero, and the 
effective number of modes for the reduced system is $12-2=10$. 
We set $m=10$ in all subsequent calculations.

For the first set of resolved modes, the reduced system includes
all the modes that are linearly unstable. For our case there are
three positive and three negative wavenumber modes that are
linearly unstable, namely the modes for $k=-3,-2,-1,1,2,3.$ 
The rest $10-6=4$ modes in the resolved set are linearly stable modes,
namely the modes with $k=-5,-4,4,5.$

\begin{figure}
\centering
\subfigure[]{\epsfig{file=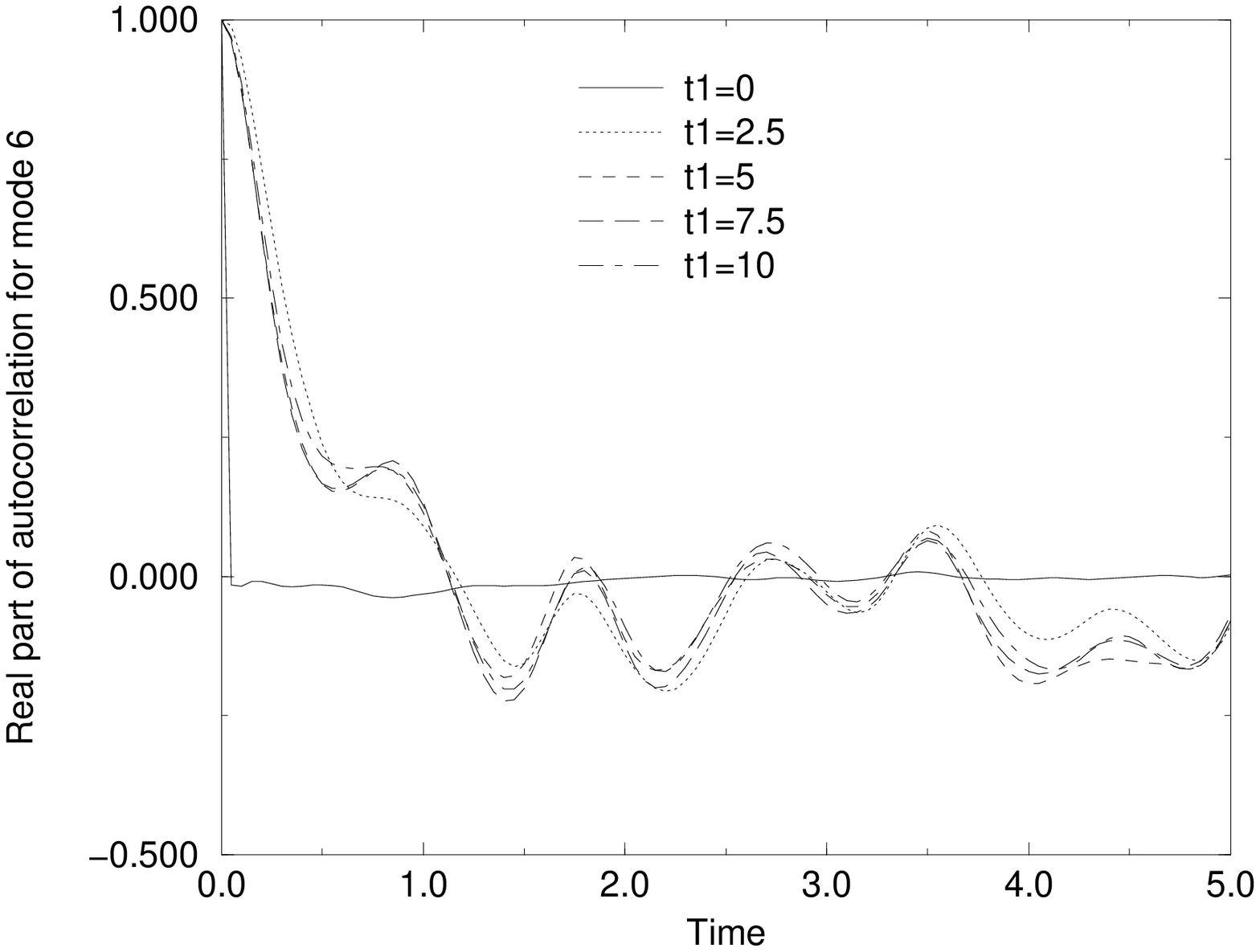,height=1.7in}}
\quad
\subfigure[]{\epsfig{file=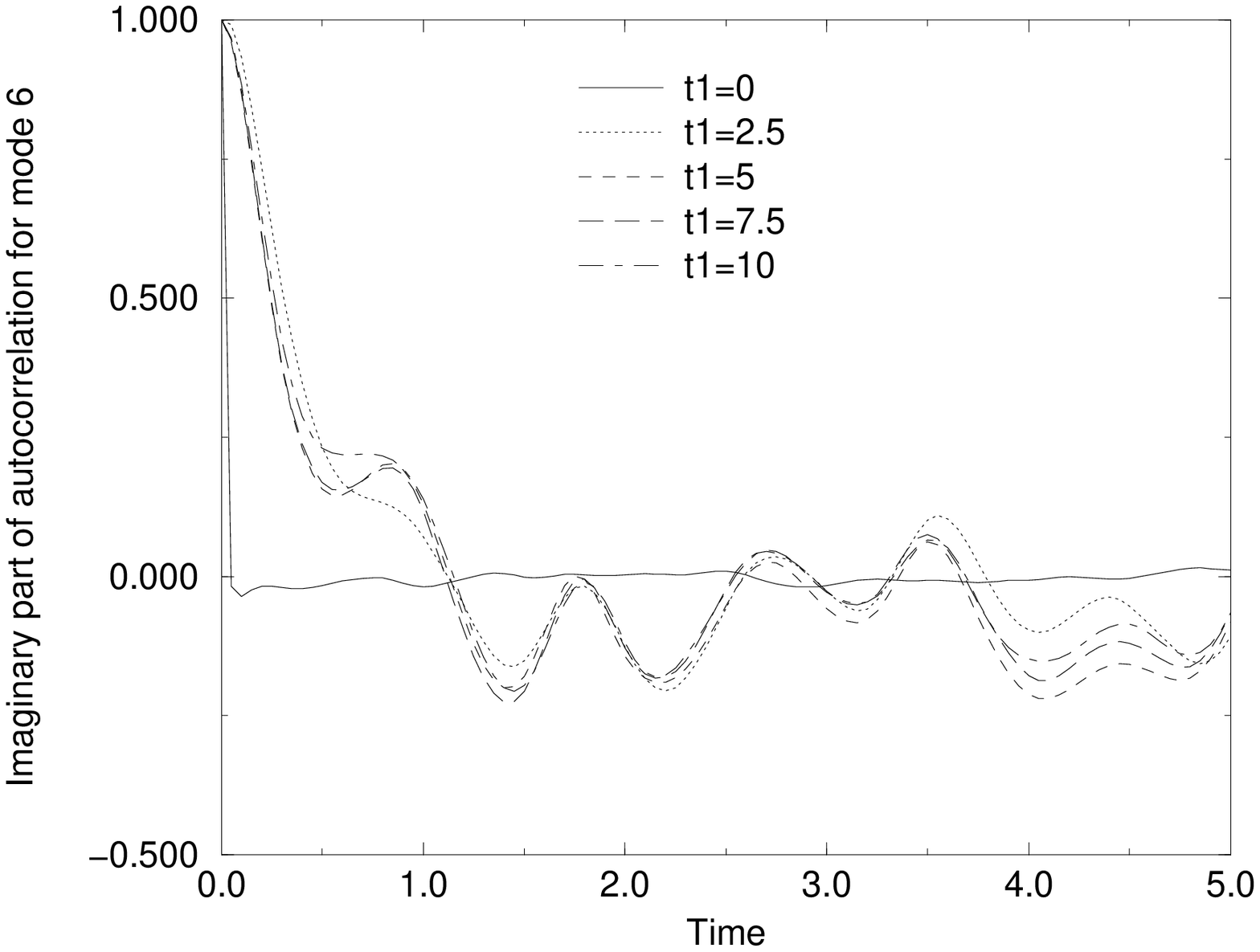,height=1.7in}}
\caption{Autocorrelation for the unresolved
mode with wavenumber 6 for different initial times. a) Real part,
b) Imaginary part.}
\label{fig:res31}
\end{figure}

\begin{figure}
\centering
\subfigure[]{\epsfig{file=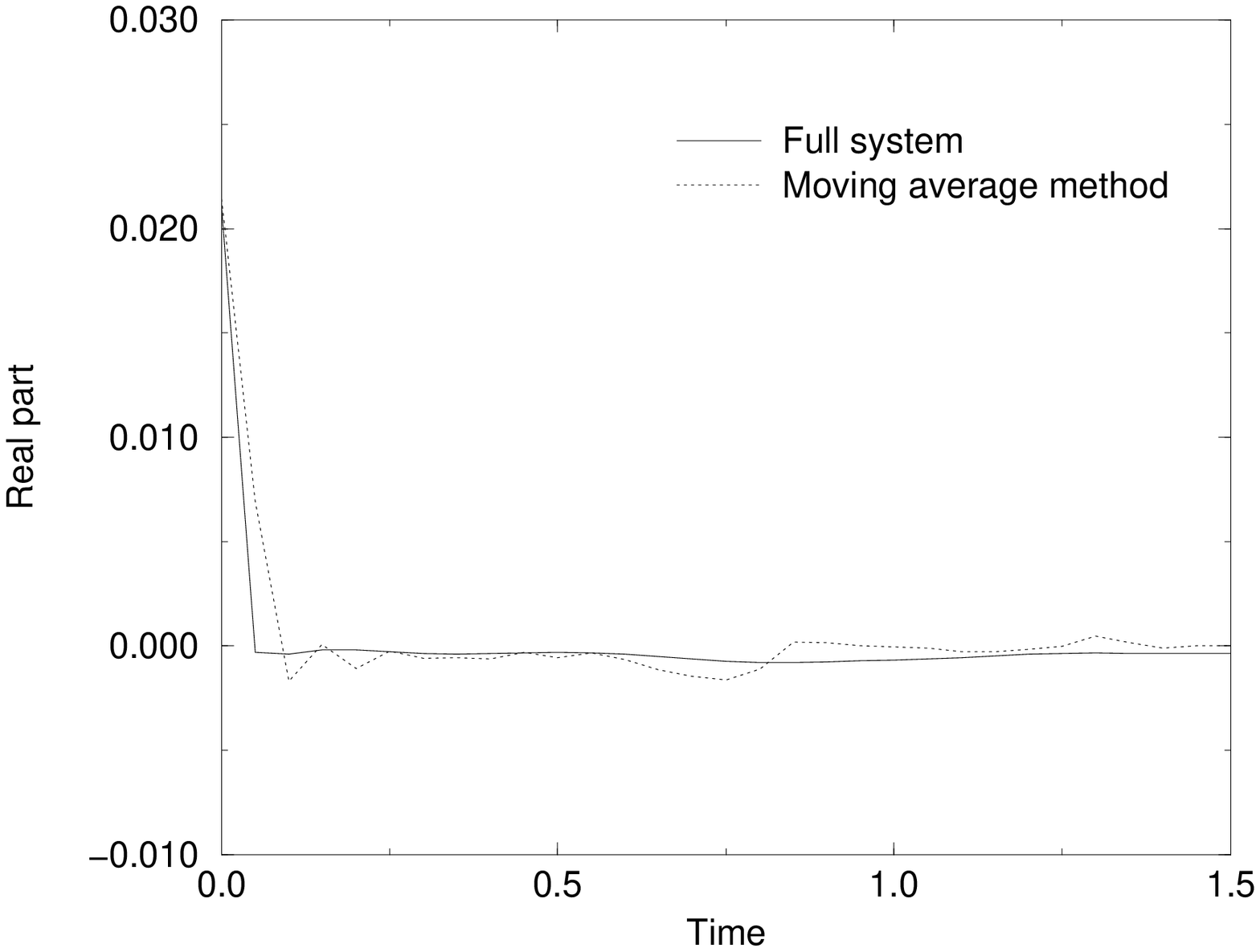,height=1.7in}}
\quad
\subfigure[]{\epsfig{file=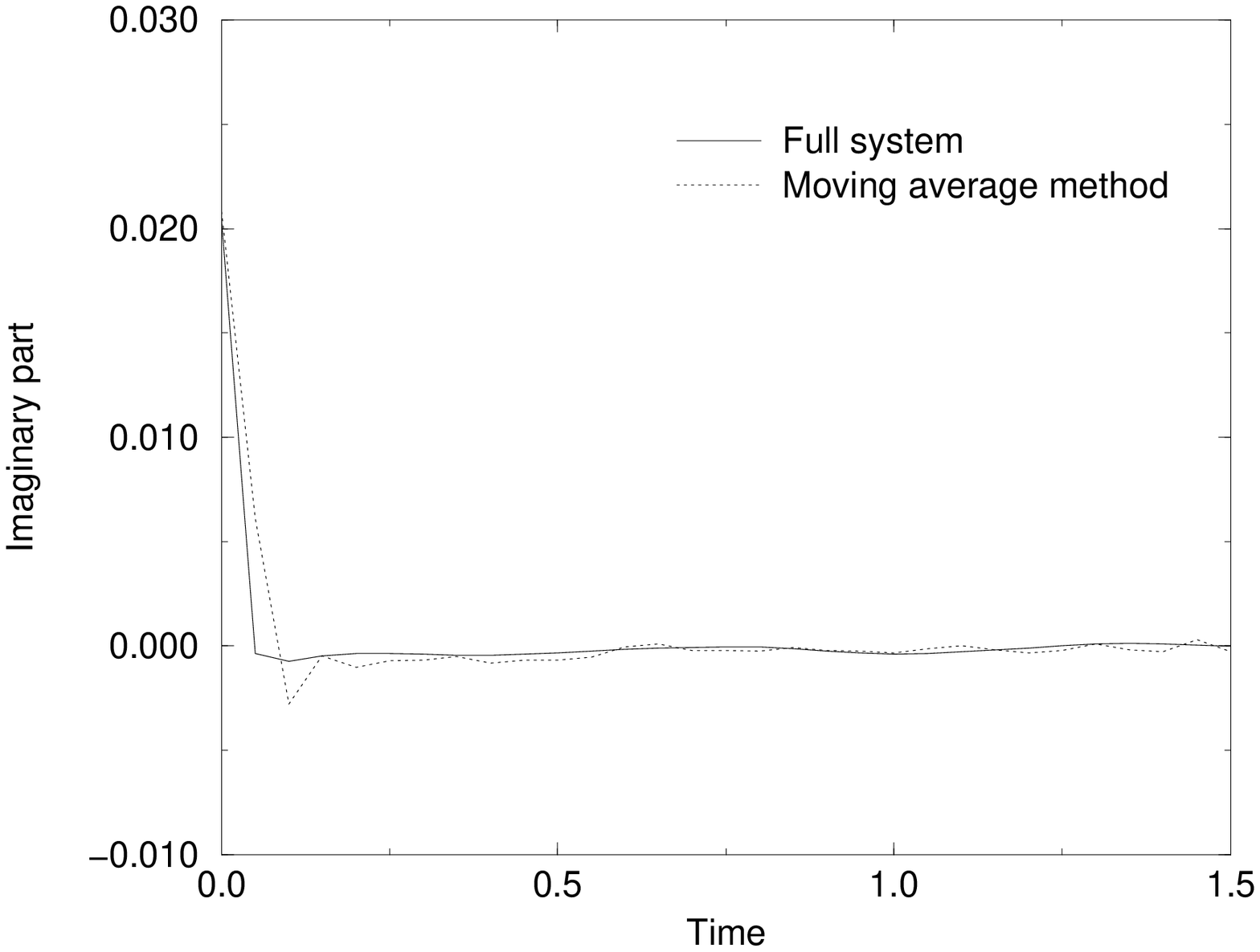,height=1.7in}}
\caption{Comparison of the autocorrelation for the unresolved
mode with wavenumber 6 as computed from the full system and 
the moving average method. a) Real part,
b) Imaginary part.}
\label{fig:res3}
\end{figure}

We proceed with the presentation of the calculations for the 
quantities needed for the evaluation of the noise term
$A_j({\fu}({\fu}_0,t)).$ The noise term
involves the unknown quantities $u_{>2}(t).$ Because we
are using a density that is not invariant, the quantities $u_{>2}(t)$ 
are nonstationary stochastic processes. Fig.(\ref{fig:res31}) shows 
the autocorrelation for the unresolved mode with wavenumber 6
$$E[(u_6(t_2)-E[u_6(t_2)])(u_6(t_1)-E[u_6(t_1)])^{*}]$$ 
for different initial times $t_1.$ As explained
before, we approximate the quantities $u_{>2}(t)$ as stationary 
stochastic processes with mean and autocorrelation determined with respect
to the diagonal Gaussian density. We use the moving average method to 
sample these stochastic processes. Each of the unknown Fourier
modes is split into its real and imaginary parts which 
are approximated independently. Fig.(\ref{fig:res3}) shows
the autocorrelation of the real and imaginary parts for the 
unresolved mode with wavenumber 6 as computed by the full system 
and the moving average method. The moving average method estimates
for the autocorrelations of the unresolved modes were produced by averaging 
over 10000 sample paths.

\begin{figure}
\centering
\subfigure[]{\epsfig{file=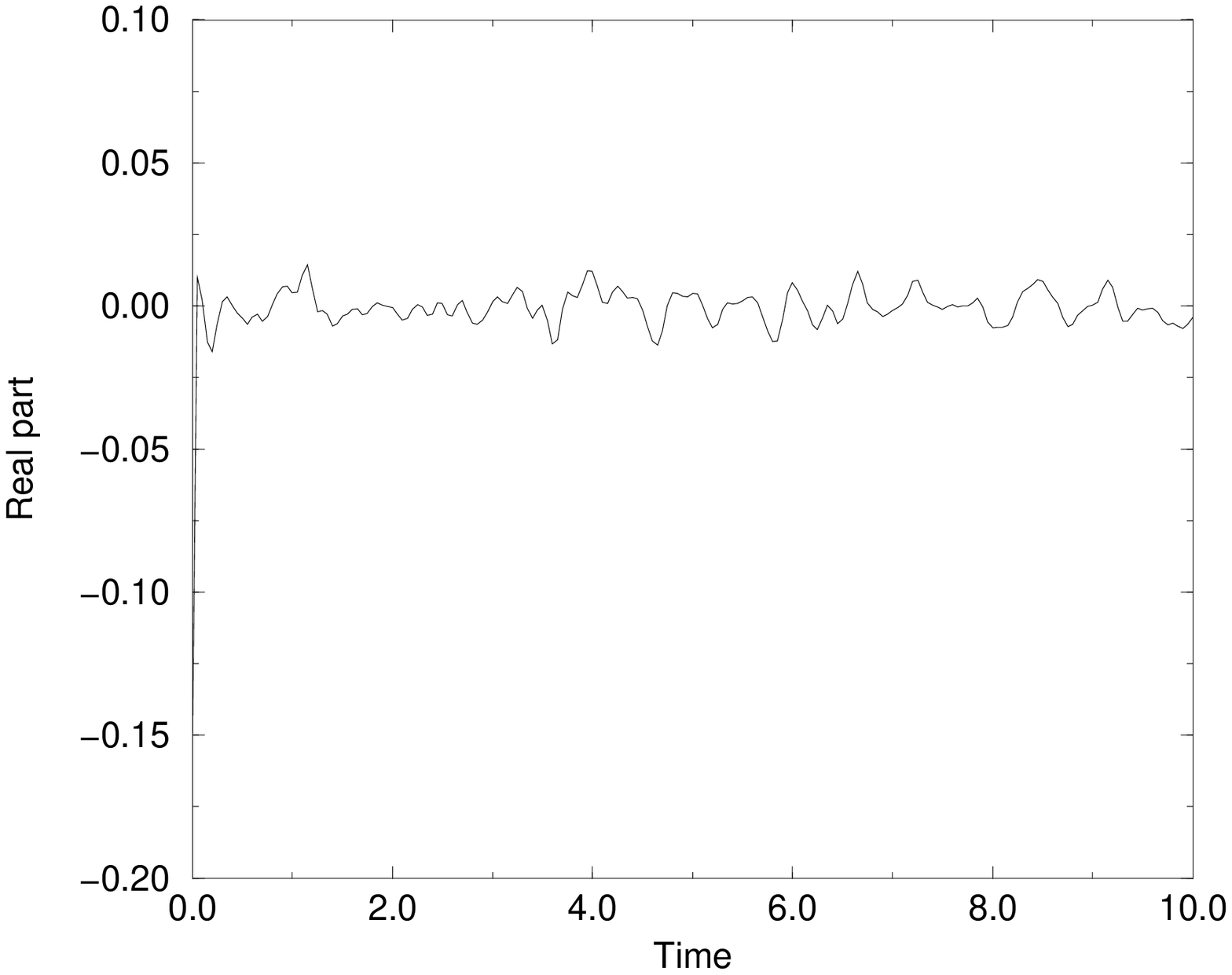,height=1.7in}}
\quad
\subfigure[]{\epsfig{file=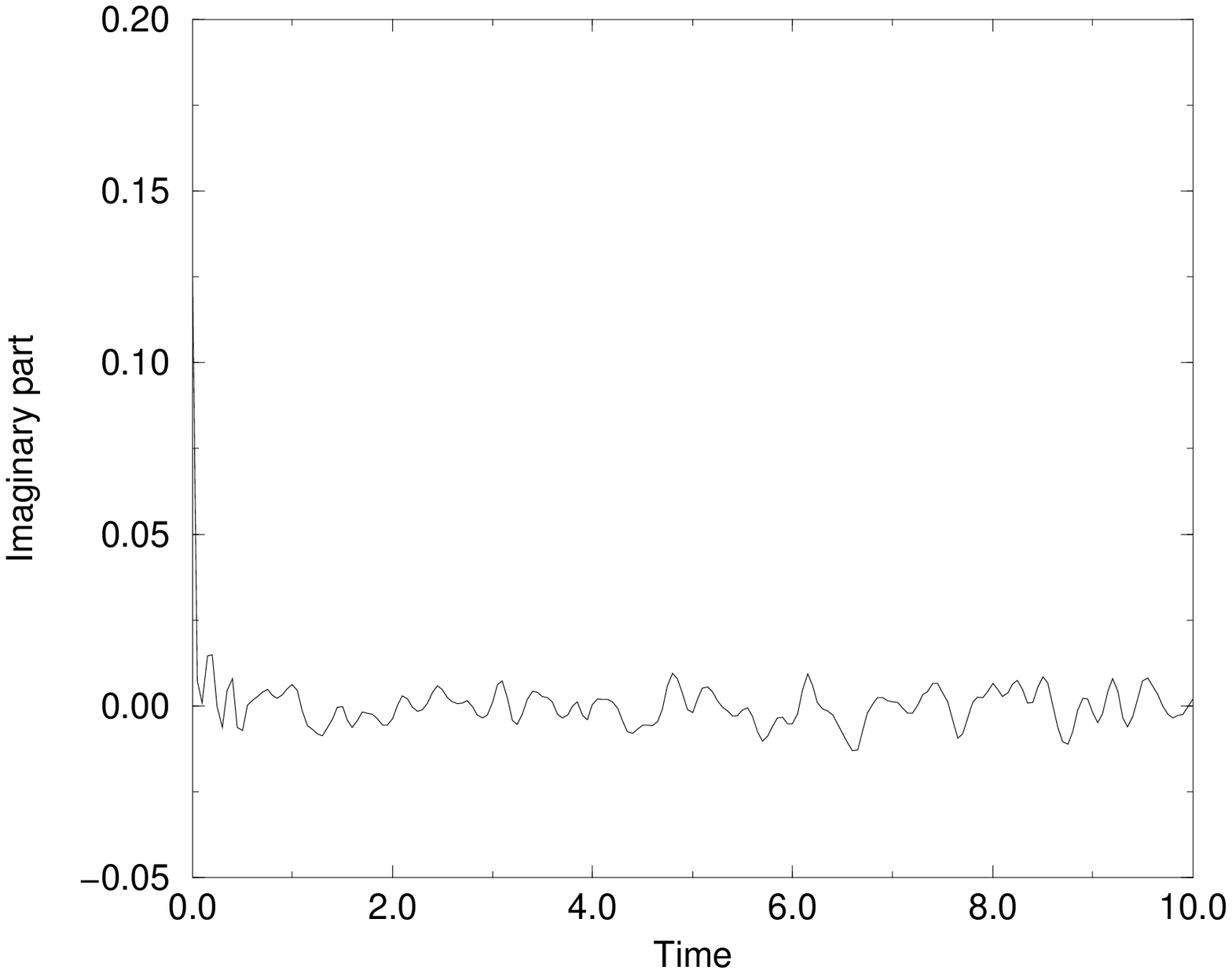,height=1.7in}}
\caption{Projection coefficient of the memory term for the
equation for the resolved mode 1 on the resolved mode 4 for
the first set of variables.
a) Real part, b) Imaginary part.}
\label{fig:res4}
\end{figure}

After the properties of the noise term, we have to calculate the 
properties of the memory term. The quantities needed for the
memory term evaluation were computed averaging over 10000 samples.
For the case of the linear projection, we have to compute the
term (recall (\ref{kspsm4}))
$$\int_0^t e^{(t-s)L} \sum_{i,k=1}^{m}b_{ik}^{-1}
(LQA_j({\fu}({\fu}_0,s)),u_{0i})u_{0k}ds,$$
so we need to compute the projection coefficients
$\sum_{i,k=1}^{m}b_{ik}^{-1}(LQA_j({\fu}({\fu}_0,s)),u_{0i}).$
Fig.(\ref{fig:res4}) shows the evolution of the 
projection coefficient  of the memory term for the equation for the resolved 
mode  with wavenumber 1  on the resolved mode 4 
$\sum_{i}^{m}b_{i4}^{-1}(LQA_1({\fu}({\fu}_0,s)),u_{0i})$ (from now on 
we omit the word wavenumber for aesthetic reasons).

Inspection of the behavior of all the projection coefficients
needed in the memory term determination reveals that
the coefficients decay fast in time. The fast decay of the quantities 
$(LQe^{sL}QLu_j,u_k)$, i.e. of the projection
coefficients in the short-memory approximation,
can only guarantee the fast decay of the quantities
$(LQe^{sQL}QLu_j,u_k)$ of the exact Mori-Zwanzig
equation for short times. For large times, the quantities
$(LQe^{sQL}QLu_j,u_k)$ can behave very differently
from $(LQe^{sL}QLu_j,u_k).$ Thus, we have to truncate the interval of
integration for the integral term from $[0,t]$ to $[t-t_0,t]$ (for $t >t_0$). 
The time $t_0$ is a short time over which the error of using 
$(LQe^{sL}QLu_j,u_k)$ instead of $(LQe^{sQL}QLu_j,u_k)$ is not large.  
Of course, we do not know the value of $t_0$ beforehand, since
we do not know the quantities $(LQe^{sQL}QLu_j,u_k).$  However, 
motivated by the numerically determined form of the projection coefficients 
$(LQe^{sL}QLu_j,u_k)$, we truncate the interval of integration for the memory
term from $[0,t]$ to $[t-1,t]$ (for $t >1$), since most projection coefficients 
decay significantly in a time interval of length 1.

The determination of the properties of the noise and the memory
allows us to solve the short-memory approximation equations for
the resolved modes. Although the full system is stiff, since it
includes many stable modes with very large decay rates,
the stable modes contained in the reduced system have small
decay rates. This allows us to use an explicit solver for the reduced
system. For our numerical experiments we use the Runge-Kutta
4th order method to solve the random integrodifferential equations
of the short-memory approximation. For the calculation of the 
integral memory term we used two different methods, namely
the trapezoidal rule and Simpson's rule. Most experiments
were done twice, once with each method.

The quantities used in the noise 
term have very fast decaying autocorrelations, which 
result in a rough (in time) noise. This, in turn, can create problems
with the numerical stability of the explicit solver. However, the noise
has small magnitude and does not appear to create problems in the
actual implementation.

The system of random integrodifferential equations is solved
over and over for different realizations of the noise. The initial
condition used for the resolved modes with wavenumbers 
$-3,-2,-1,1,2,3$ is 
$u_1=u_{-1}=u_2=u_{-2}=u_3=u_{-3}=1$, while the rest
of the resolved modes are set initially equal to zero.
The conditional expectations of the resolved modes conditioned
on their initial values are computed by averaging over the
realizations of the noise.

The truth is computed as follows. We sample the
density over and over keeping the resolved modes fixed to the values just
described. For each sample, we evolve (\ref{ksp}) for $k=-11,\ldots,11.$ 
Finally, we average over the samples and obtain the conditional expectations of 
the resolved modes.

The Galerkin approximation consists of setting the unresolved modes 
equal to zero for all times, i.e. solving (\ref{ksp}) for $k=-5,\ldots,5.$ 
For the Galerkin approximation there is no need to solve the reduced system 
repeatedly, because for each realization the unresolved modes are set equal to 
zero.

\begin{figure}
\centering
\subfigure[]{\epsfig{file=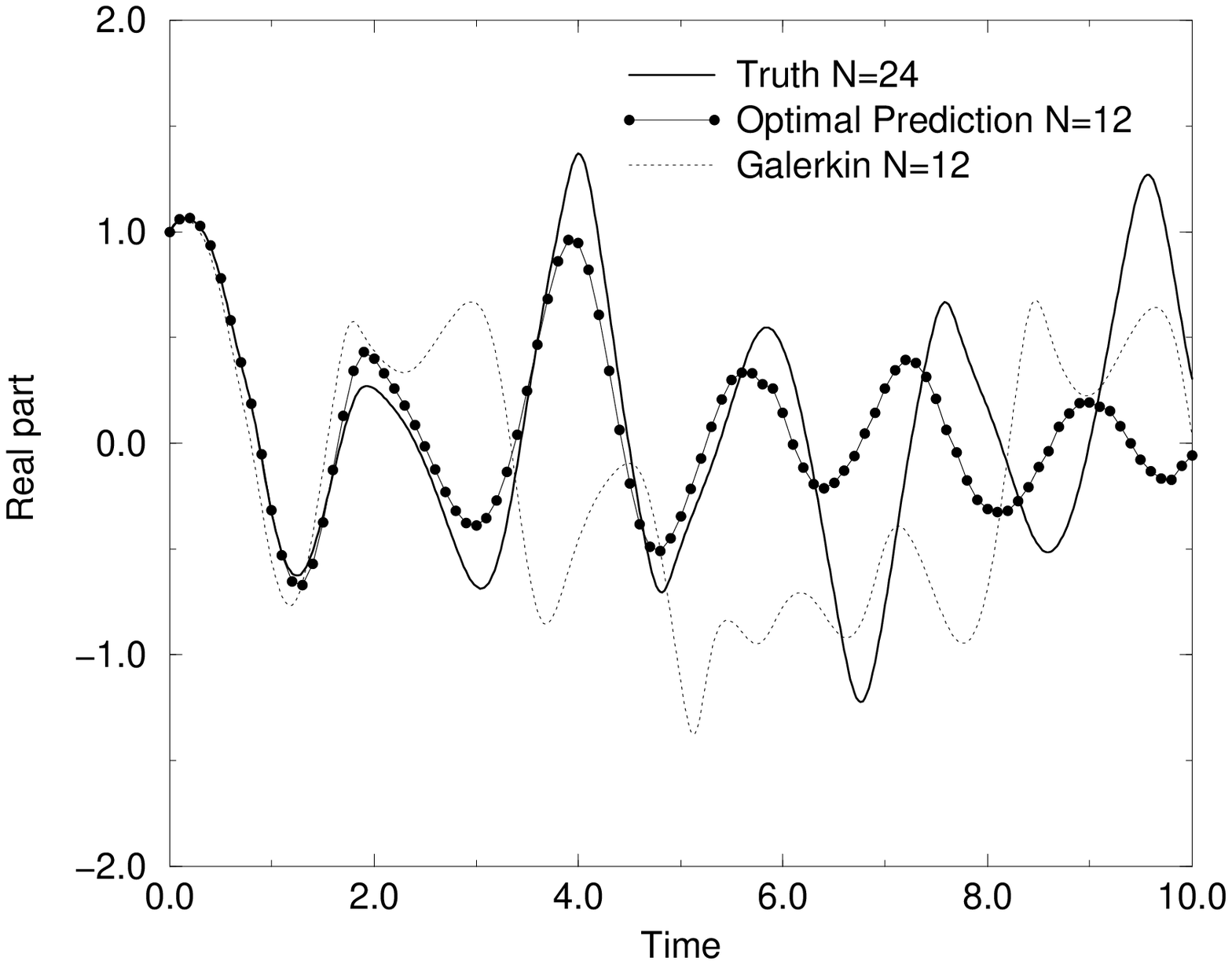,height=1.7in}}
\quad
\subfigure[]{\epsfig{file=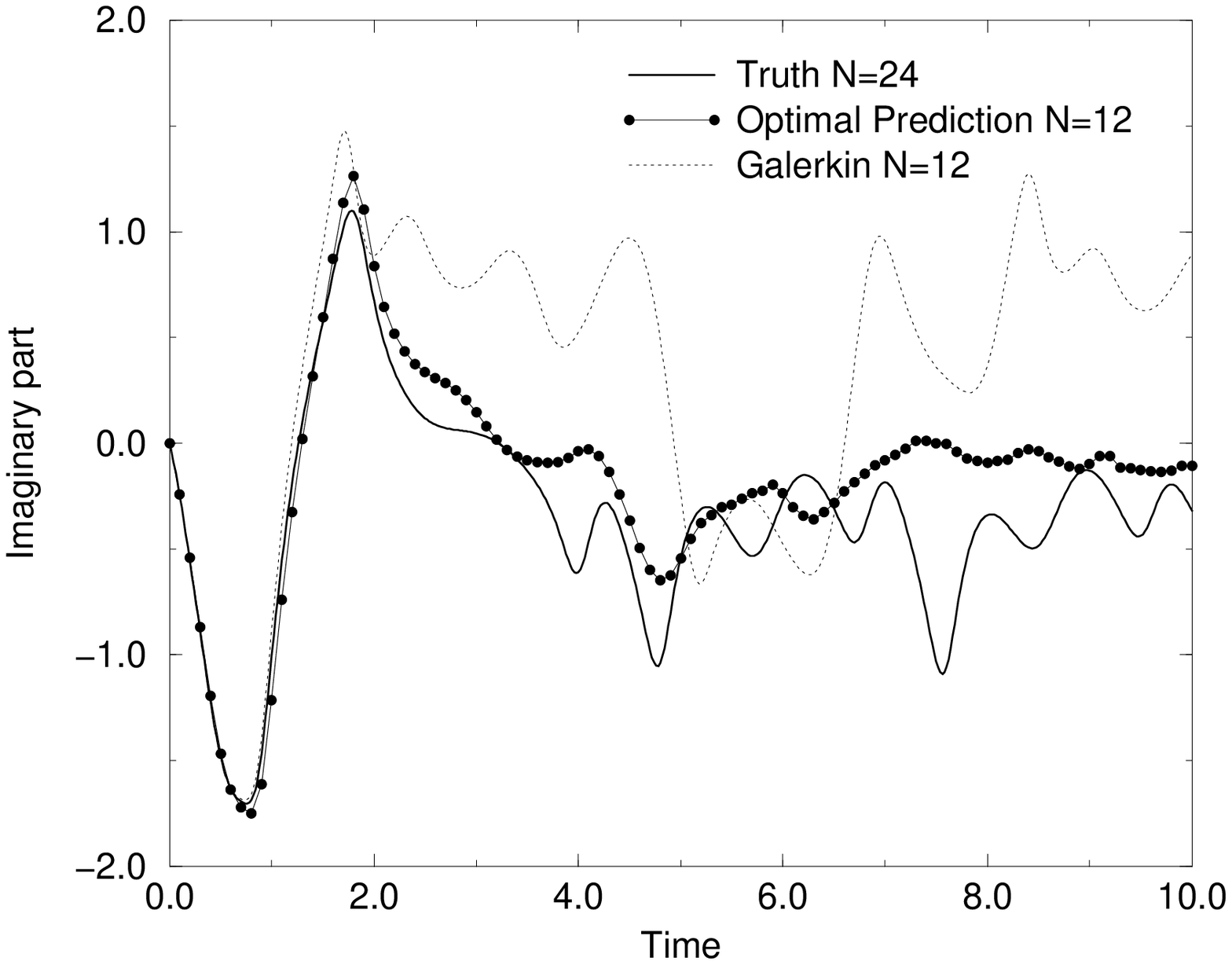,height=1.7in}}
\caption{Conditional expectation evolution for the
 resolved mode 1 for the first set of resolved variables. a) Real part, 
 b) Imaginary part.}
\label{fig:res6}
\end{figure}

\begin{figure}
\centering
\subfigure[]{\epsfig{file=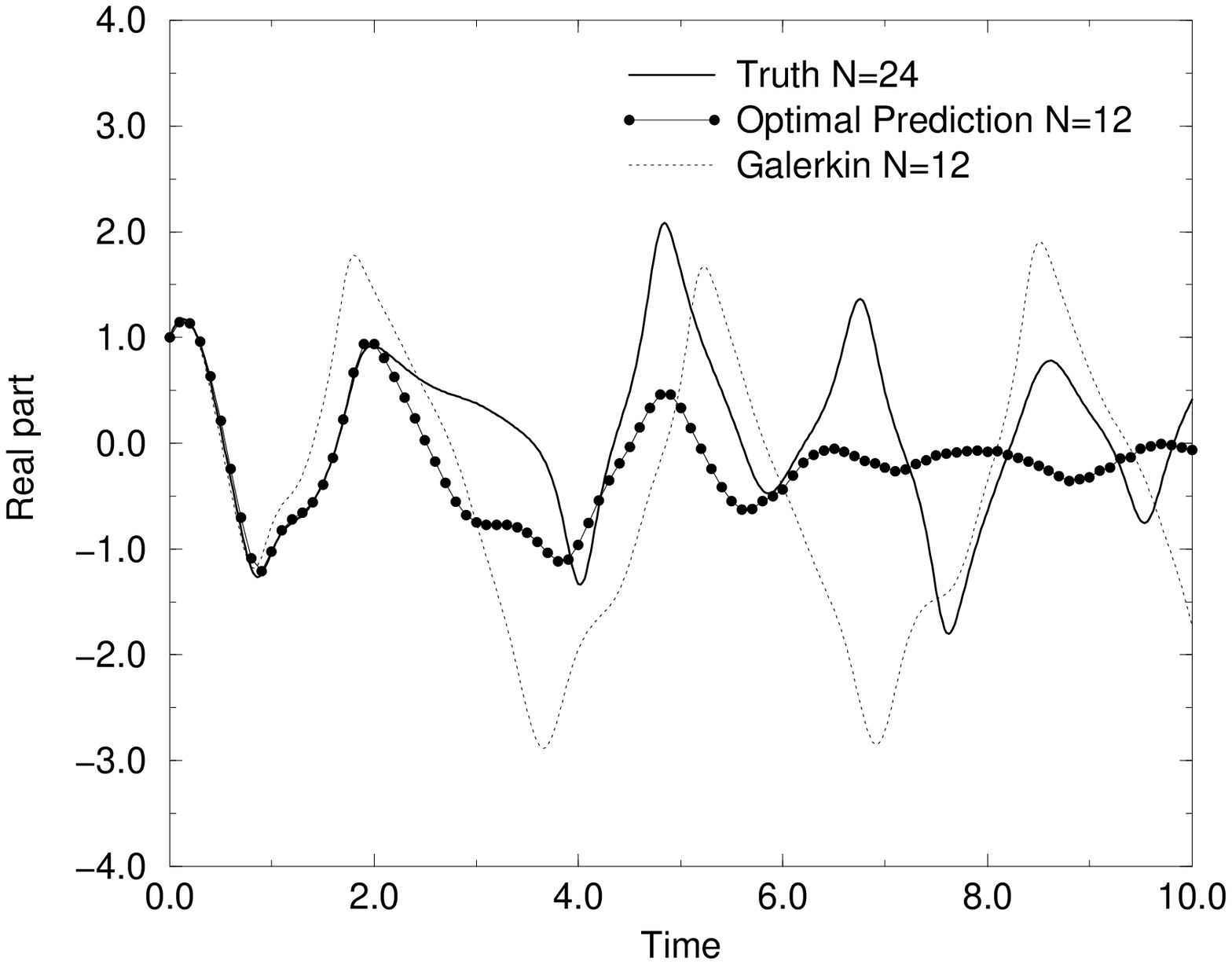,height=1.7in}}
\quad
\subfigure[]{\epsfig{file=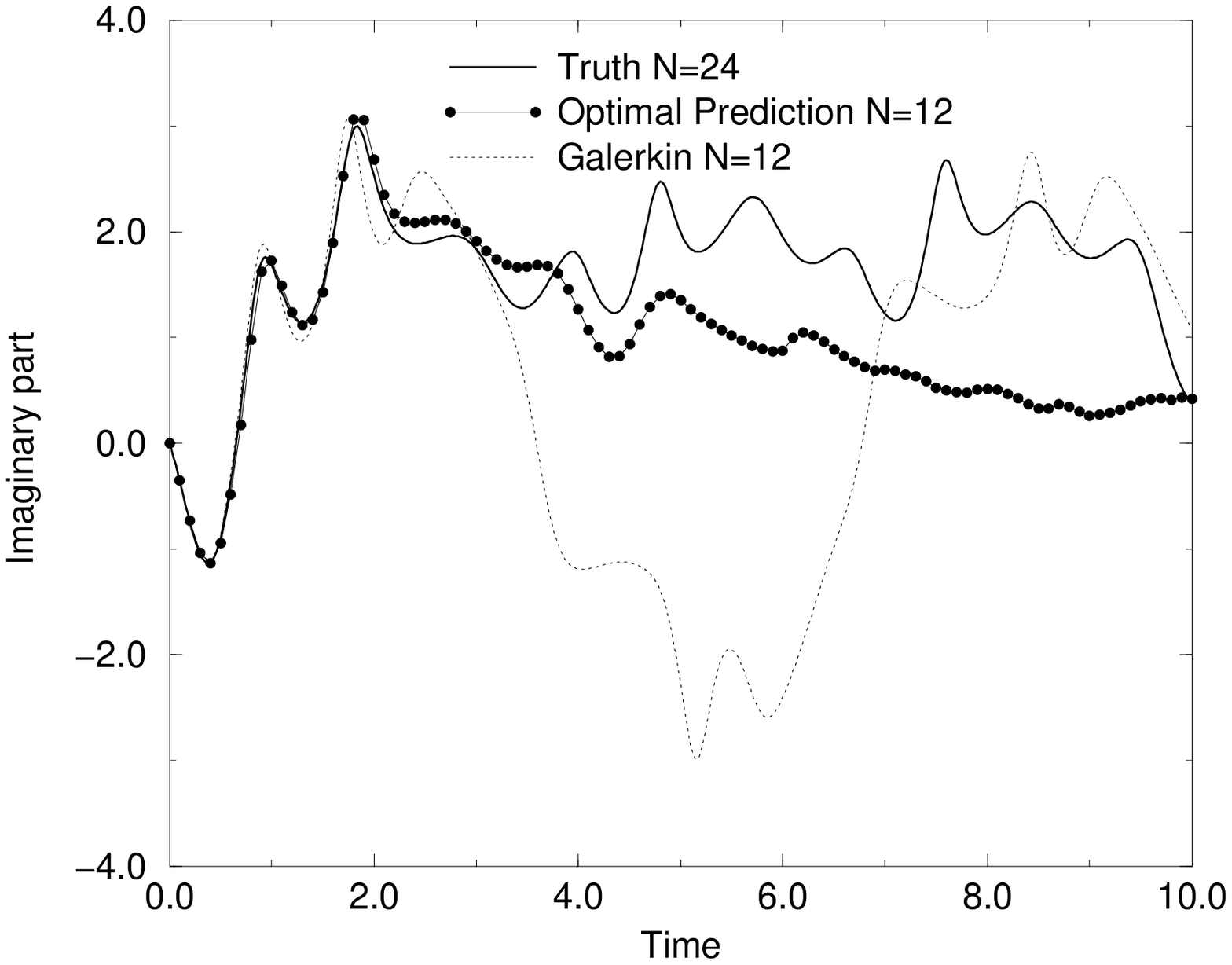,height=1.7in}}
\caption{Conditional expectation evolution for the
 resolved mode 2 for the first set of resolved variables. 
 a) Real part, b) Imaginary part.}
\label{fig:res7}
\end{figure}

\begin{figure}
\centering
\subfigure[]{\epsfig{file=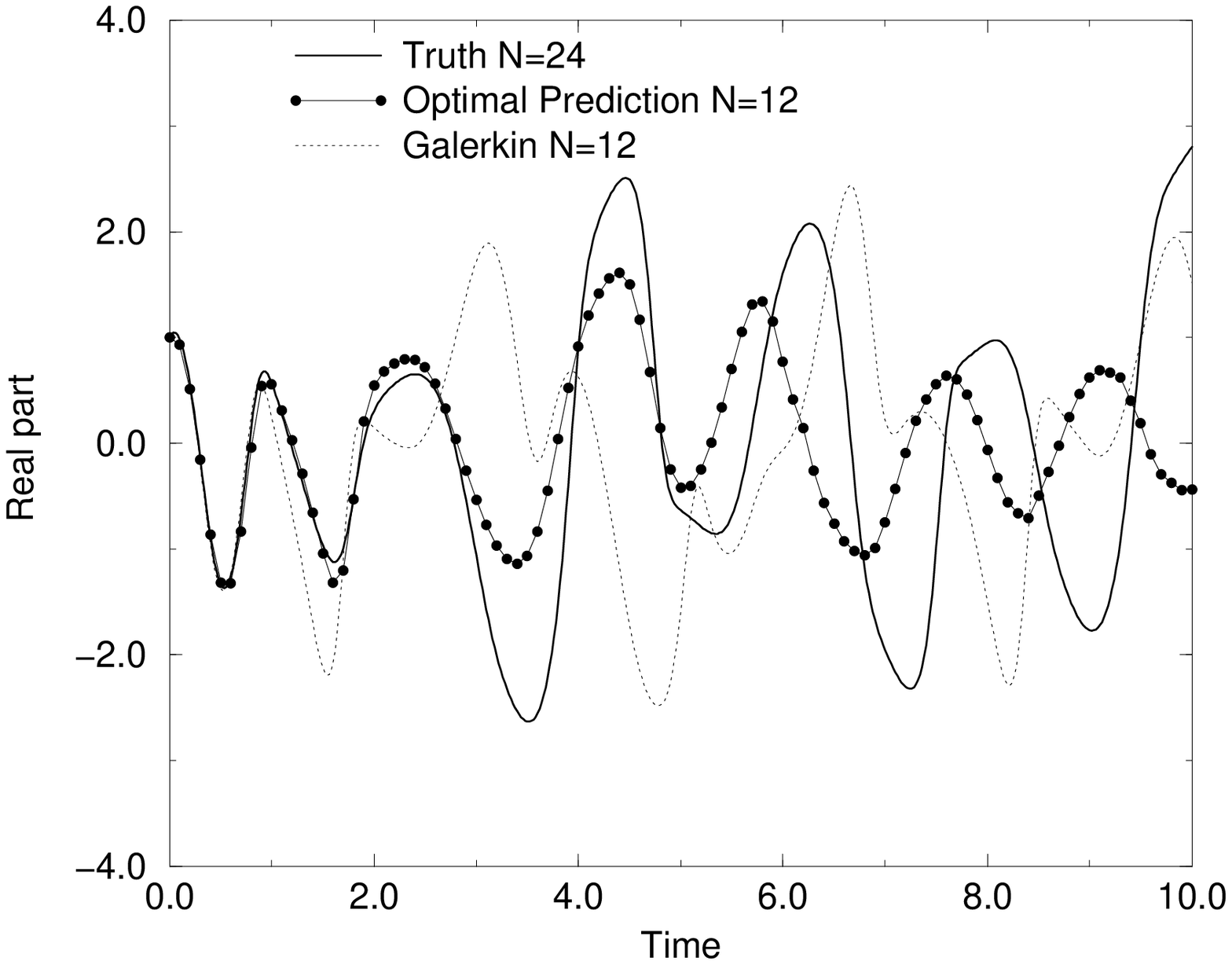,height=1.7in}}
\quad
\subfigure[]{\epsfig{file=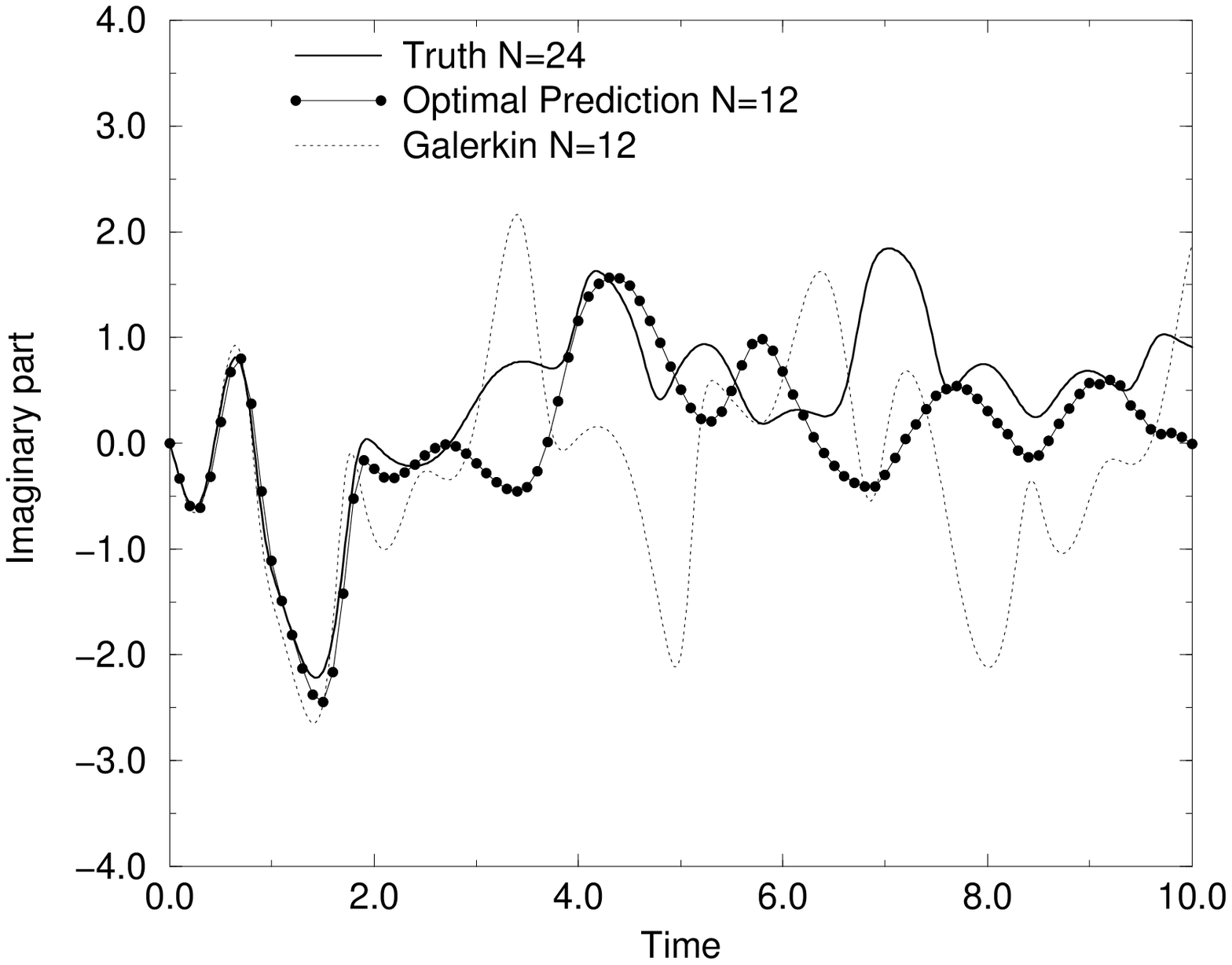,height=1.7in}}
\caption{Conditional expectation evolution for the
 resolved mode 3 for the first set of resolved variables. 
 a) Real part, b) Imaginary part.}
\label{fig:res8}
\end{figure}

Figs.(\ref{fig:res6}), (\ref{fig:res7}), 
(\ref{fig:res8}) show the real and imaginary parts of the short-memory 
approximation estimates of the
conditional expectations for the resolved modes 1,2,3
as compared to the truth and the Galerkin approximation. The
true conditional expectations were computed by averaging over
1000 samples. The short-memory estimates were computed
by averaging over 1000 realizations of the noise. The results for the
truth and the optimal prediction estimates are converged, even with
a 1000 samples, because, the unresolved modes have small magnitudes.
However, there is no "slaving" \cite{fo1} of the unresolved modes to the 
resolved ones. The randomness in the initial conditions of the unresolved modes 
does affect the conditional expectations of the resolved ones. However, in the 
time interval that we examine, the results of the randomness of the initial 
conditions have not yet fully manifested themselves, and thus 1000 samples 
(or 1000 realizations of the noise) are enough for 
convergence of the results. For longer times, more samples (or realizations
of the noise) are needed for convergence of the truth and the optimal
prediction estimates.

As can be seen, 
there is good agreement between
the short-memory approximation estimates for the conditional
expectations of the resolved modes and the true conditional
expectations of these modes. 
Also, the short-memory approximation results 
constitute a considerable improvement over the results for
the Galerkin approximation. The good agreement of the 
short-memory estimates for the conditional expectations
of the resolved modes, indicates that the fast decay of
the projection coefficients $(LQe^{sQL}QLu_j,u_k)$
is not restricted to short times only. On the contrary, they
must continue to decay fast for larger times and this allows 
the short-memory approximation to perform well.

\begin{figure}
\centering
\subfigure[]{\epsfig{file=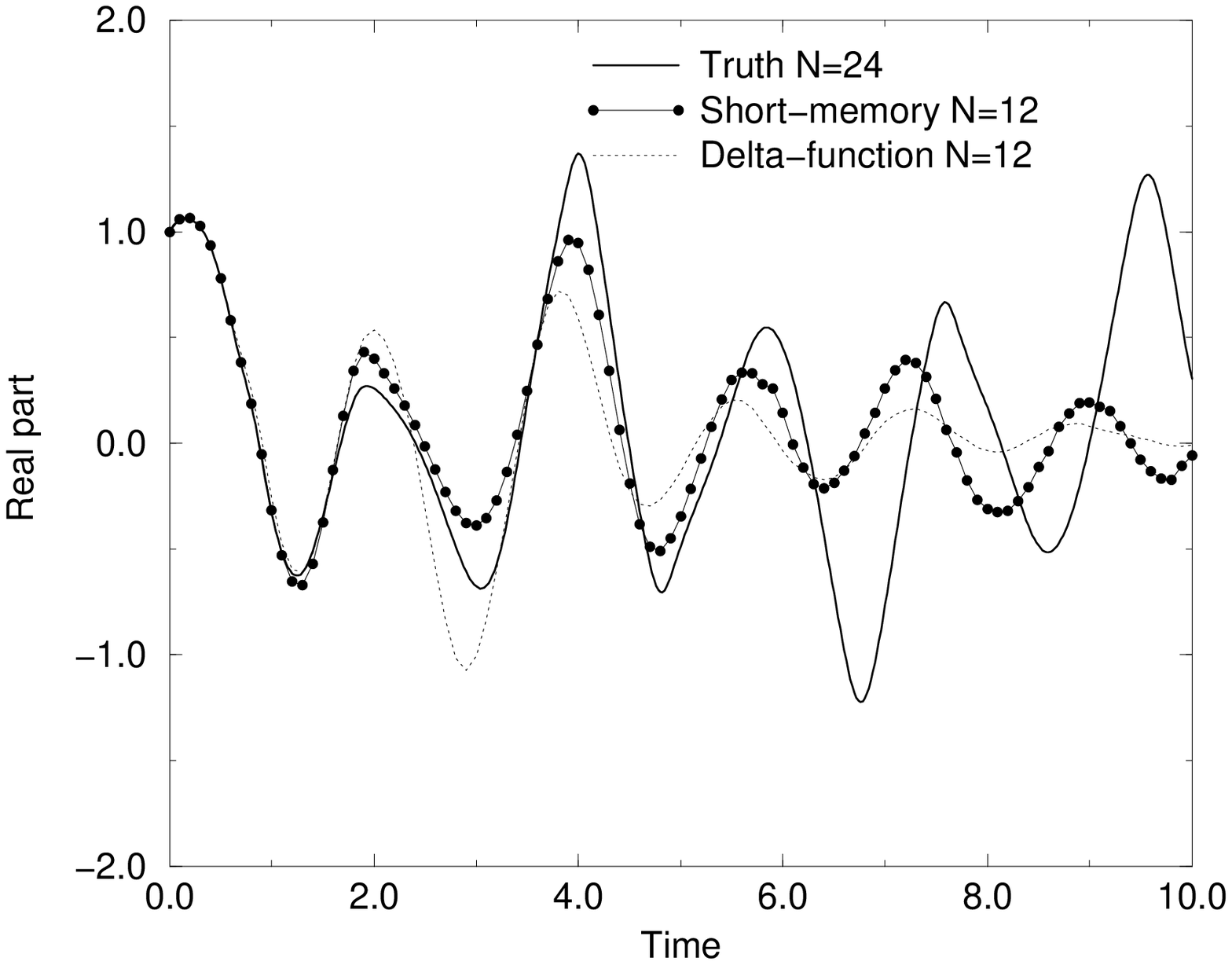,height=1.7in}}
\quad
\subfigure[]{\epsfig{file=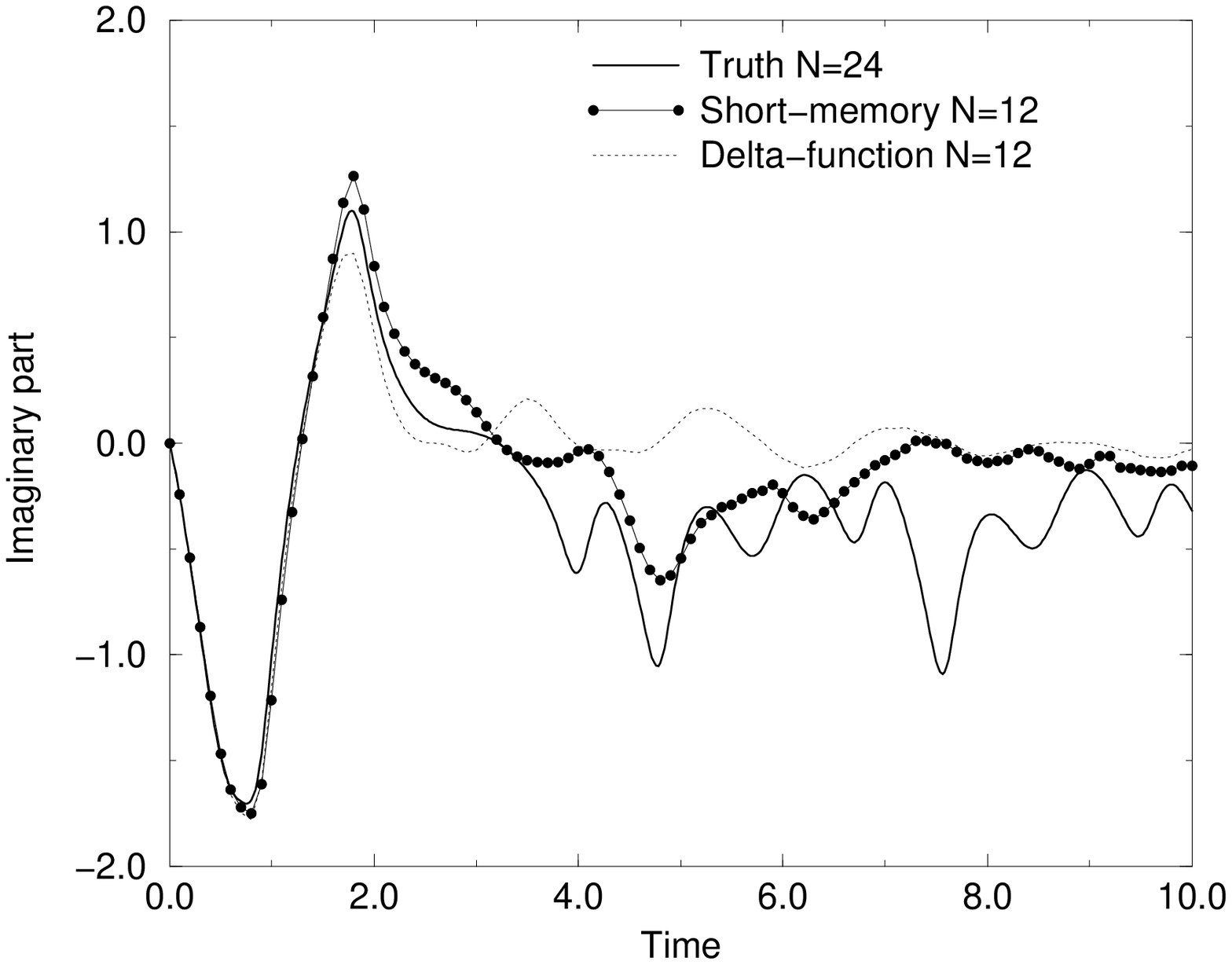,height=1.7in}}
\caption{Conditional expectation evolution for the
resolved mode 1 for the first set of resolved variables. 
Delta-function vs short-memory. 
a) Real part, b) Imaginary part.}
\label{fig:res91}
\end{figure}

\begin{figure}
\centering
\subfigure[]{\epsfig{file=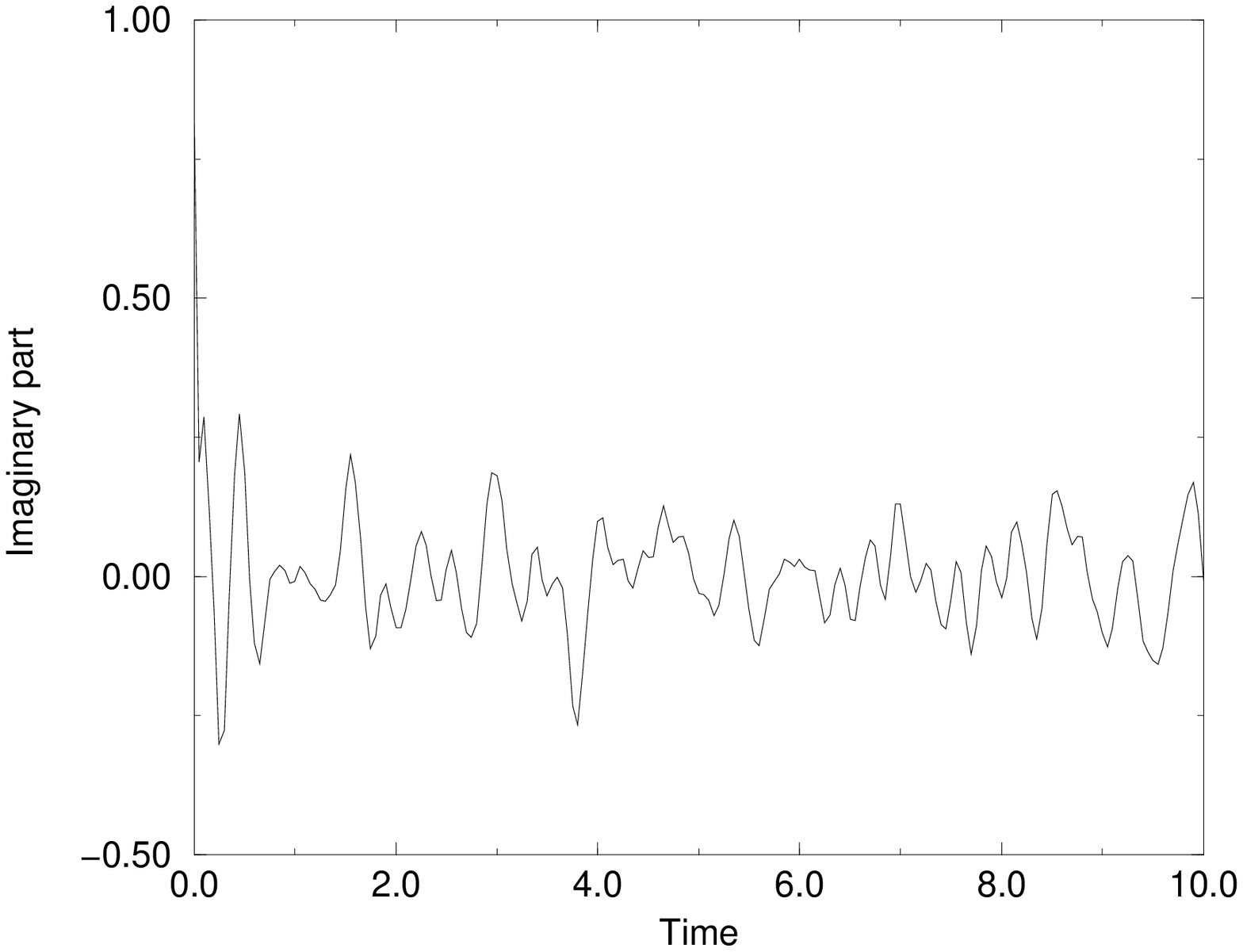,height=1.7in}}
\caption{Imaginary part of the projection coefficient of the memory term for the
equation for the resolved mode 5 on the resolved mode 1 
for the first set of resolved variables.}
\label{fig:res10}
\end{figure}

Motivated by the fact that all the projection coefficients
$(LQe^{sL}QLu_j,u_k)$
used in the short-memory approximation decay fast,
we can inquire about the validity of an even more
drastic approximation which replaces
the projection coefficients by delta-functions multiplied
by the integrals of the quantities $(LQe^{sL}QLu_j,u_k).$
Of course, as before, these integrals must be computed
over a truncated interval. The advantage of such an
approximation is that the reduced system is a system
of random ordinary differential equations and not
integrodifferential equations. This results in a significant
reduction of the computational time needed to calculate
the estimates of the conditional expectations of the resolved
modes.  Fig.(\ref{fig:res91}) show the estimate of the
conditional expectation for the resolved mode 1
using the delta-function approximation, compared to the results of the 
short-memory approximation and the truth. The delta-function
approximation can predict accurately the evolution of the conditional 
expectation for shorter times than the short-memory approximation. 
The inadequacy of the delta-function
approximation  to predict the evolution of the conditional expectations for
longer times can be explained by a more careful inspection of the projection 
coefficients $(LQe^{sL}QLu_j,u_k).$ Although these coefficients decay fast, 
several of them do not retain the same sign during this decay 
(see e.g. Fig.(\ref{fig:res10})). In fact, they oscillate rapidly with 
non-negligible amplitudes. Integration of these oscillations results
in cancellations that lead to loss of important information about the
short-time memory. As a result, the delta-function approximation loses
accuracy. However, the accuracy exhibited by the delta-function approximation
is impressive, if we consider that it results in a system of differential (and
not integrodifferential) 
differential equations whose numerical integration is very fast. For  our
numerical experiments, the integration of the delta-function approximation
equations is 15 times faster than the integration of the short-memory
approximation equations. In fact, the time of integration, for one realization of
the noise, of the delta-function approximation equations, is comparable to the 
time of integration of the equations for the simple Galerkin approximation. 

\subsubsection{Finite-rank projection}

We conclude our presentation of results for the first set of variables with
the estimates of the conditional expectations obtained when we use the 
finite-rank projection in the memory term calculations.

\begin{figure}
\centering
\subfigure[]{\epsfig{file=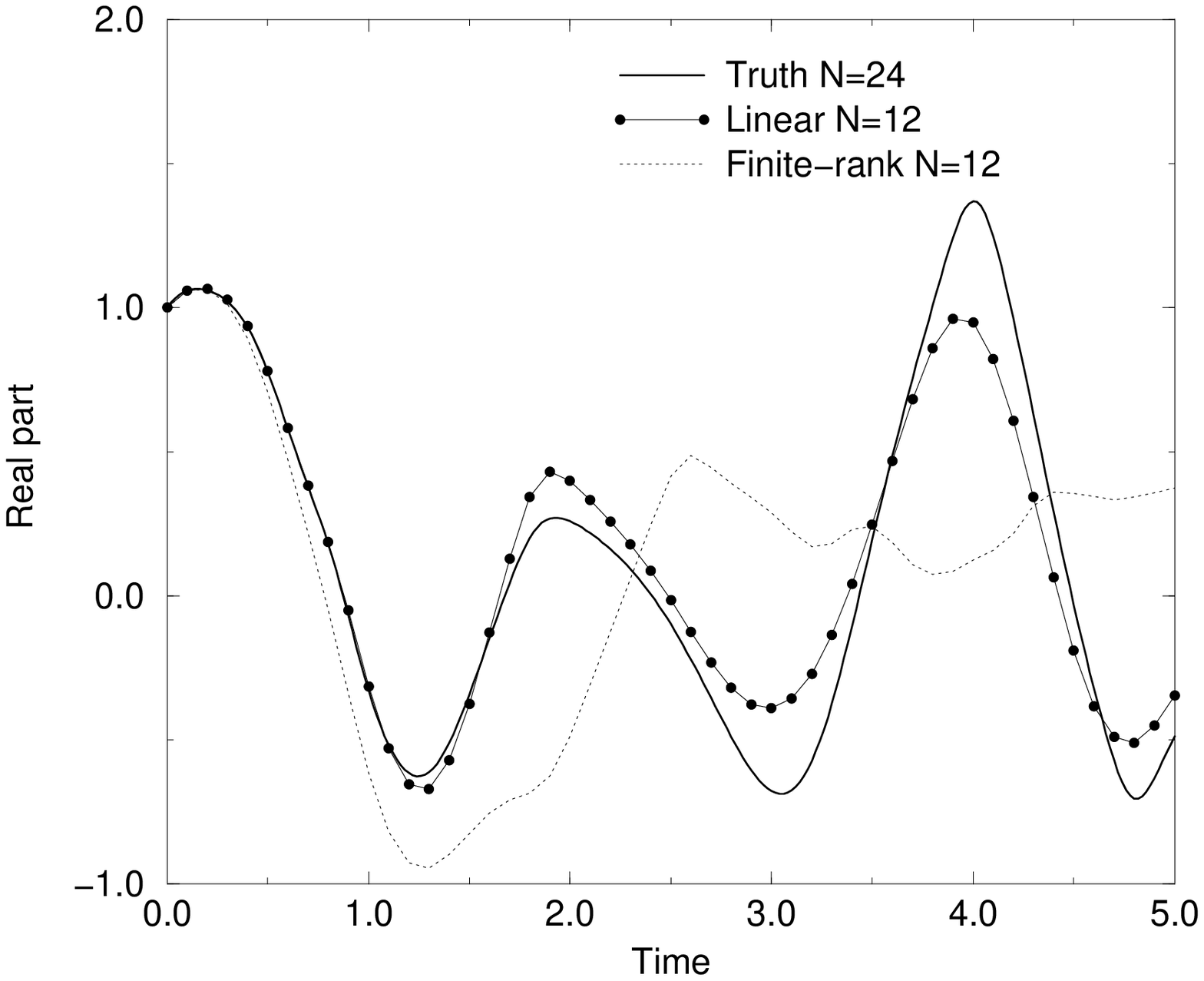,height=1.7in}}
\quad
\subfigure[]{\epsfig{file=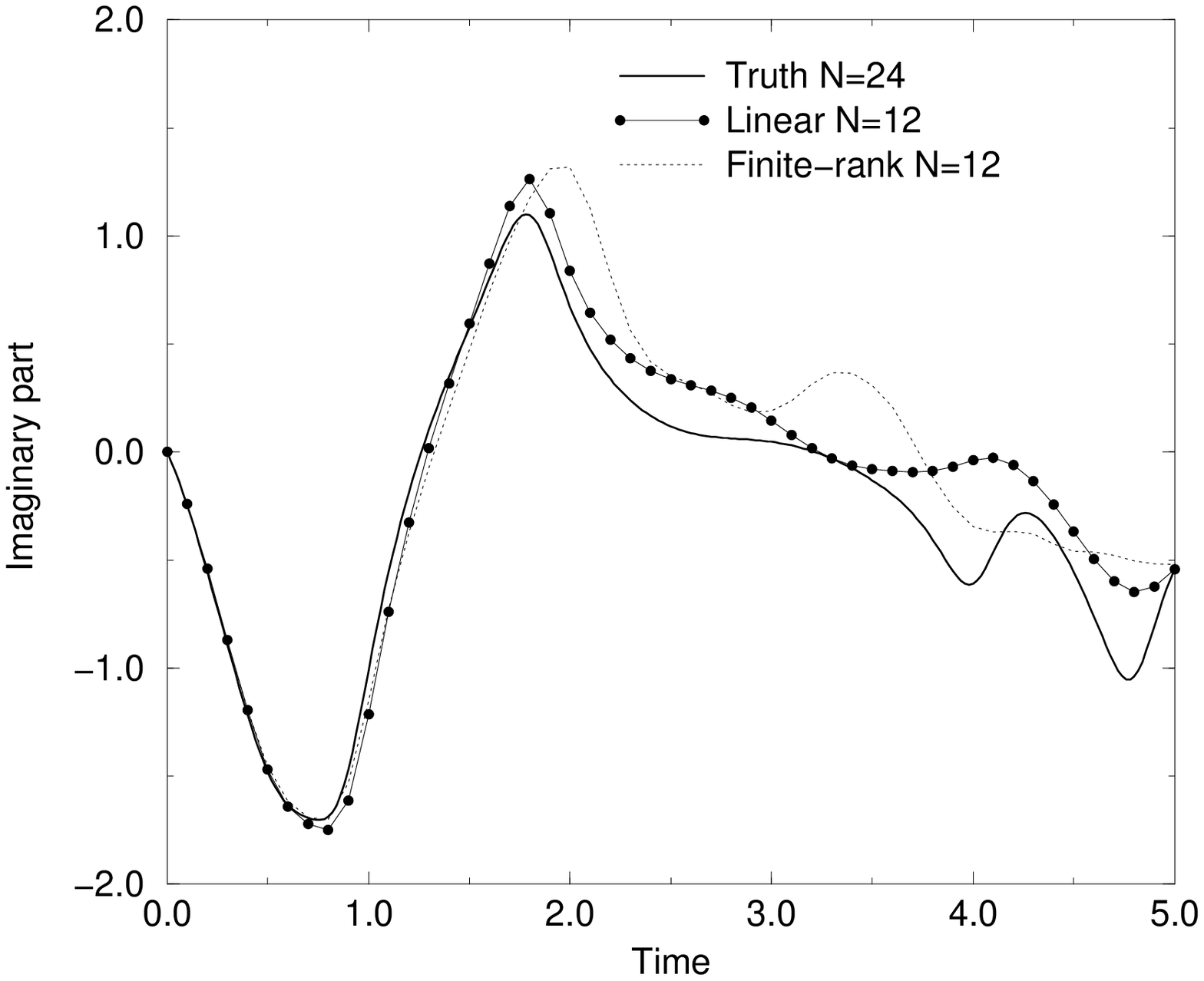,height=1.7in}}
\caption{Conditional expectation evolution for the
 resolved mode 1 for the first set of resolved variables. Linear
 vs. finite-rank projection for the memory. a) Real part, b) Imaginary 
part.}
\label{fig:res10f1}
\end{figure}

Fig.(\ref{fig:res10f1}) shows the real and imaginary parts of the 
finite-rank projection short-memory approximation estimates of the
conditional expectations for the resolved mode 1, as compared to the truth 
and the linear projection short-memory approximation 
estimates. The finite-rank projection short-memory approximation estimates 
were computed by averaging over 1000 realizations of the noise. For our 
numerical experiments we use the value $\beta_{\kappa_j}=0$ for all 
$H_{\kappa_j}$, irrespective of their order (see \ref{finite1}-\ref{finite3}). 
Note that when $\beta_{\kappa_j}=0,$ the first order polynomials in the set 
are the functions used in the linear projection (with a different  numerical 
factor coming from the orthonormality of the set functions).
The finite-rank projection estimates are computed using polynomials of order
up to 2 for the two most unstable modes. All the constraints mentioned above
result in a set of 30 orthonormal functions (we also omit the function that 
consists of polynomials of order zero in each mode, which is a constant). The
interval of integration of the memory term was truncated from $[0,t]$ to 
$[t-1,t]$. Different
truncations did not alter the results much.

\begin{figure}
\centering
\subfigure[]{\epsfig{file=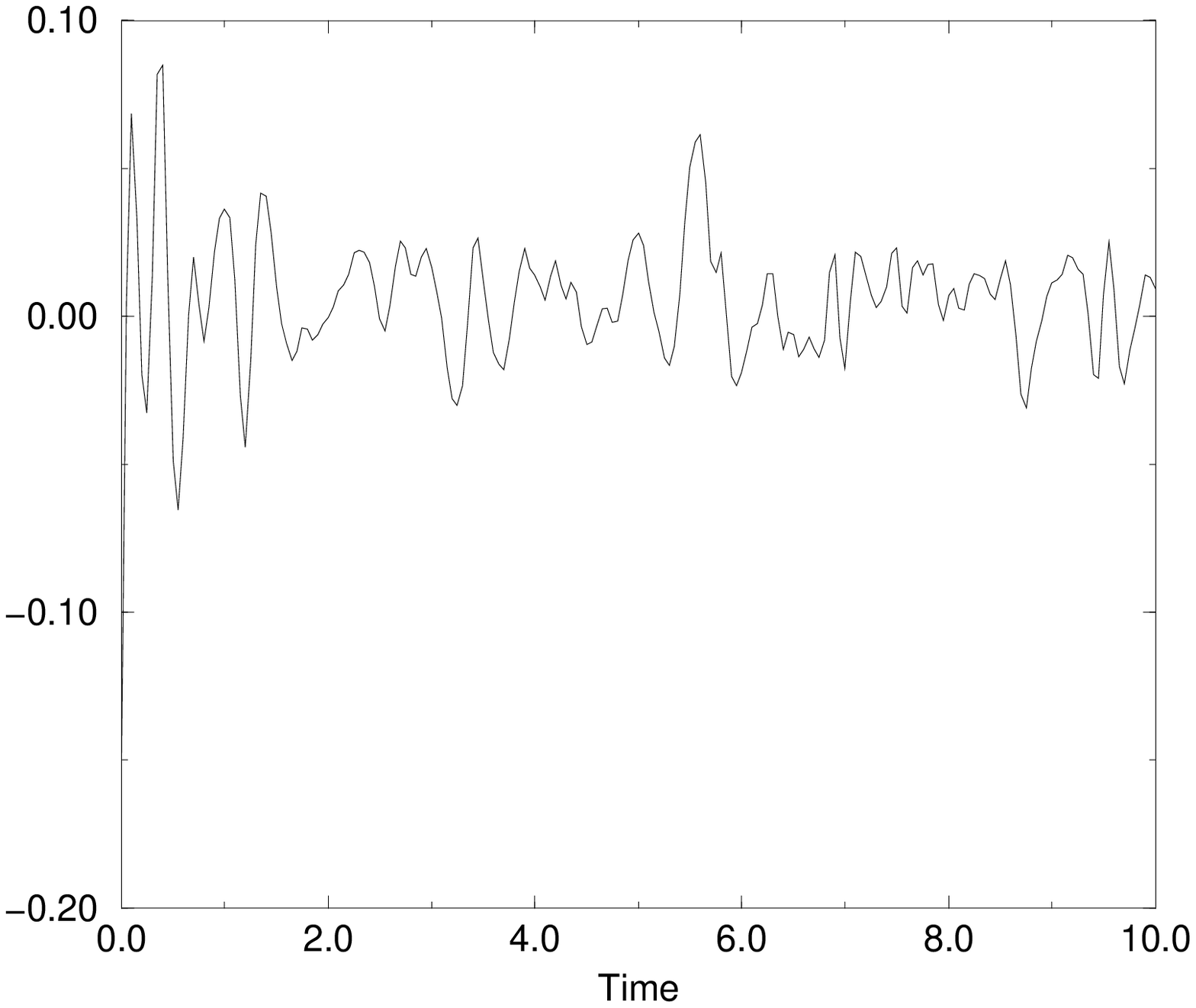,height=1.7in}}
\quad
\subfigure[]{\epsfig{file=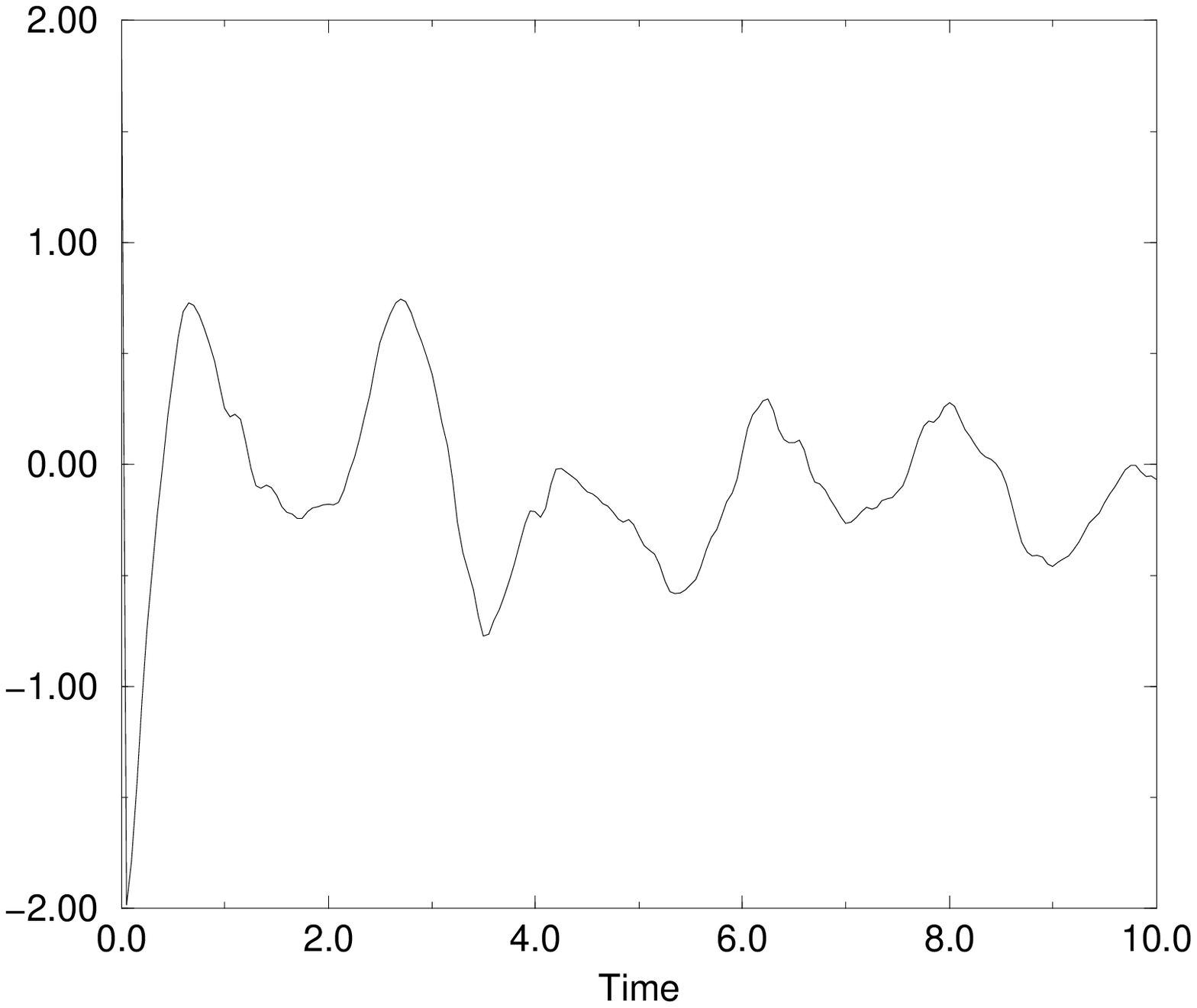,height=1.7in}}
\caption{Examples of slowly-decaying projection coefficients
of the memory term on polynomials of order higher than 1. Finite-rank
projection for the first set of resolved variables.}
\label{fig:res10f2}
\end{figure}

The finite-rank projection estimates of the conditional expectations are worse 
than the linear projection ones. This result, although surprising at first, can
be explained by examining the assumption under which the short-memory
approximation is expected to be valid for long times. The short-memory 
approximation is valid if the quantities $(LQe^{sQL}QLu_j,h_k)$ decay fast 
and stay small for long times. If they do not, then the short-memory 
approximation is valid for short times only. Also, the behavior of the quantities 
$(LQe^{sL}QLu_j,h_k)$ can only be used to infer the behavior of 
$(LQe^{sQL}QLu_j,h_k)$ for short times. Inspection of the projection 
coefficients $(LQe^{sL}QLu_j,h_k)$ for the case of the finite-rank projection 
reveals that several of these coefficients, do not decay fast. This
means that, at least for short times, the quantities $(LQe^{sQL}QLu_j,h_k)$
do not decay fast. Since the conditional expectation estimates based on the
finite-rank projection are not good for longer times, we conclude that the
quantities $(LQe^{sQL}QLu_j,h_k)$ do not start decaying fast after a short
time. This is not contradictory with the notion that the finite-rank projection
is a more accurate projection than the linear one. Indeed, the finite rank
projection is better than the linear one for the memory term of the exact 
Mori-Zwanzig equation and for the memory term of the short-memory 
approximation. 
However, being a more accurate projection, does not necessarily guarantee that 
the projection coefficients on polynomials of order higher than 1 
(recall $\beta_{\kappa_j}=0$) have to decay fast and stay small for long times. 
This is exactly what happens in our case. Several projection coefficients 
$(LQe^{sQL}QLu_j,h_k)$ on polynomials 
of order higher than 1 do not decay fast and do not stay small for long times,
thus violating the short-memory approximation assumption. In this way, the
short-memory approximation for the linear projection case that does not
include those slowly-decaying coefficients, is better than the short-memory
approximation for the finite-rank projection that does include those
slowly-decaying projection coefficients. Fig.(\ref{fig:res10f2}) shows
some examples of slowly decaying projection coefficients on polynomials of 
order larger than 1.

\subsection{Resolution of all but one unstable modes}

For the second set of resolved modes, we include in the
resolved modes all the linearly unstable modes except
for one. We choose to leave the modes $k=-1,1$ as
unresolved. There is no particular reason for this
choice. Any other unstable mode can be left as unresolved, since
it is the instability, and not the growth rate, of an unresolved mode that 
determines the accuracy of the short-memory approximation. In fact, 
for our value of the viscosity, mode 1 has the smallest (linear) growth rate 
among the unstable modes. With this choice, the set of resolved modes
contains the unstable modes $k=-3,-2,2,3$ and the stable
modes $k=-6,-5,-4,4,5,6.$ As in the case of the first set of resolved variables, 
the full system has $n=22$ modes and the reduced system $m=10$ modes.
For the second set of resolved variables, we only conducted 
numerical experiments where the linear projection was used
for the memory term.

\begin{figure}
\centering
\subfigure[]{\epsfig{file=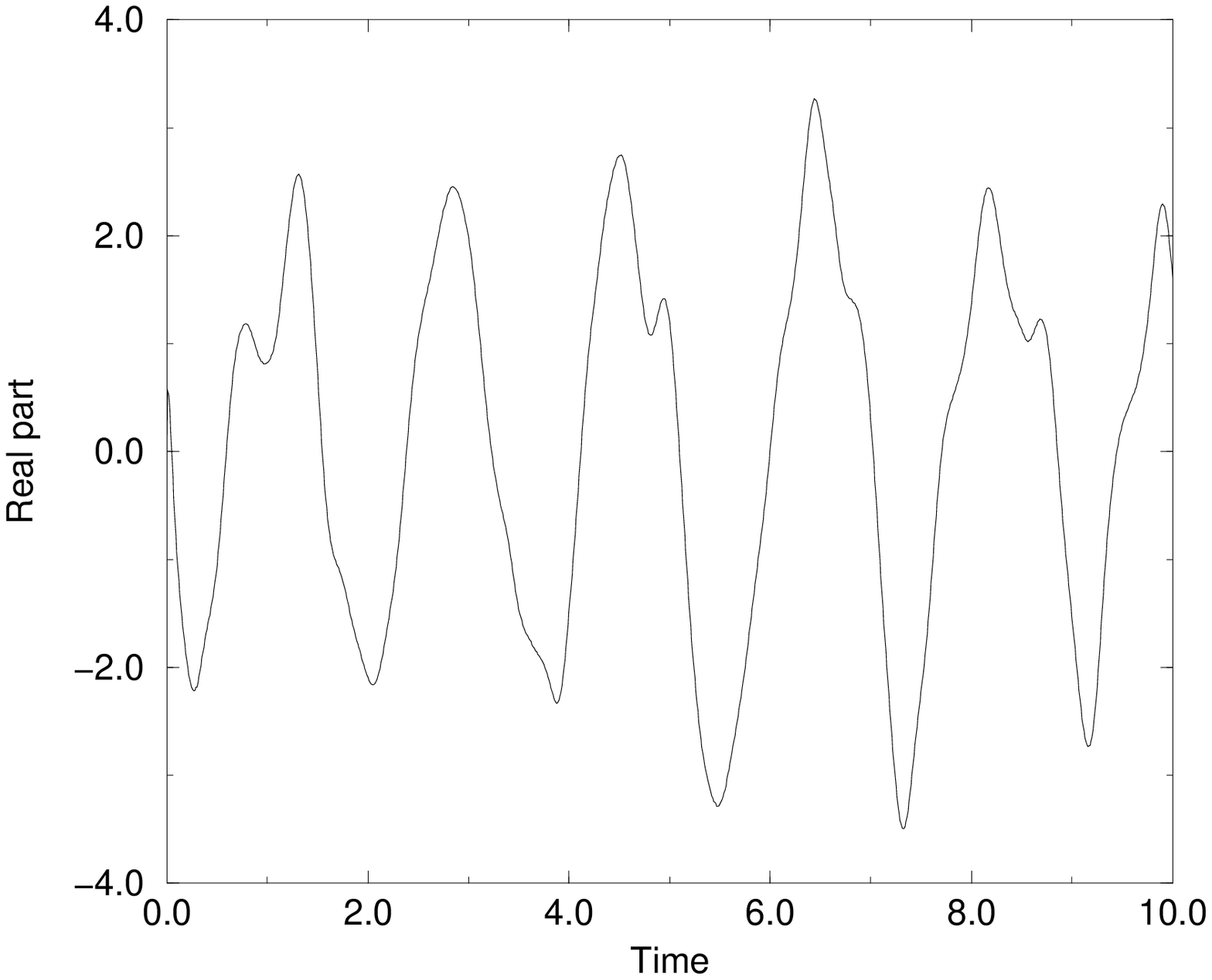,height=1.7in}}
\quad
\subfigure[]{\epsfig{file=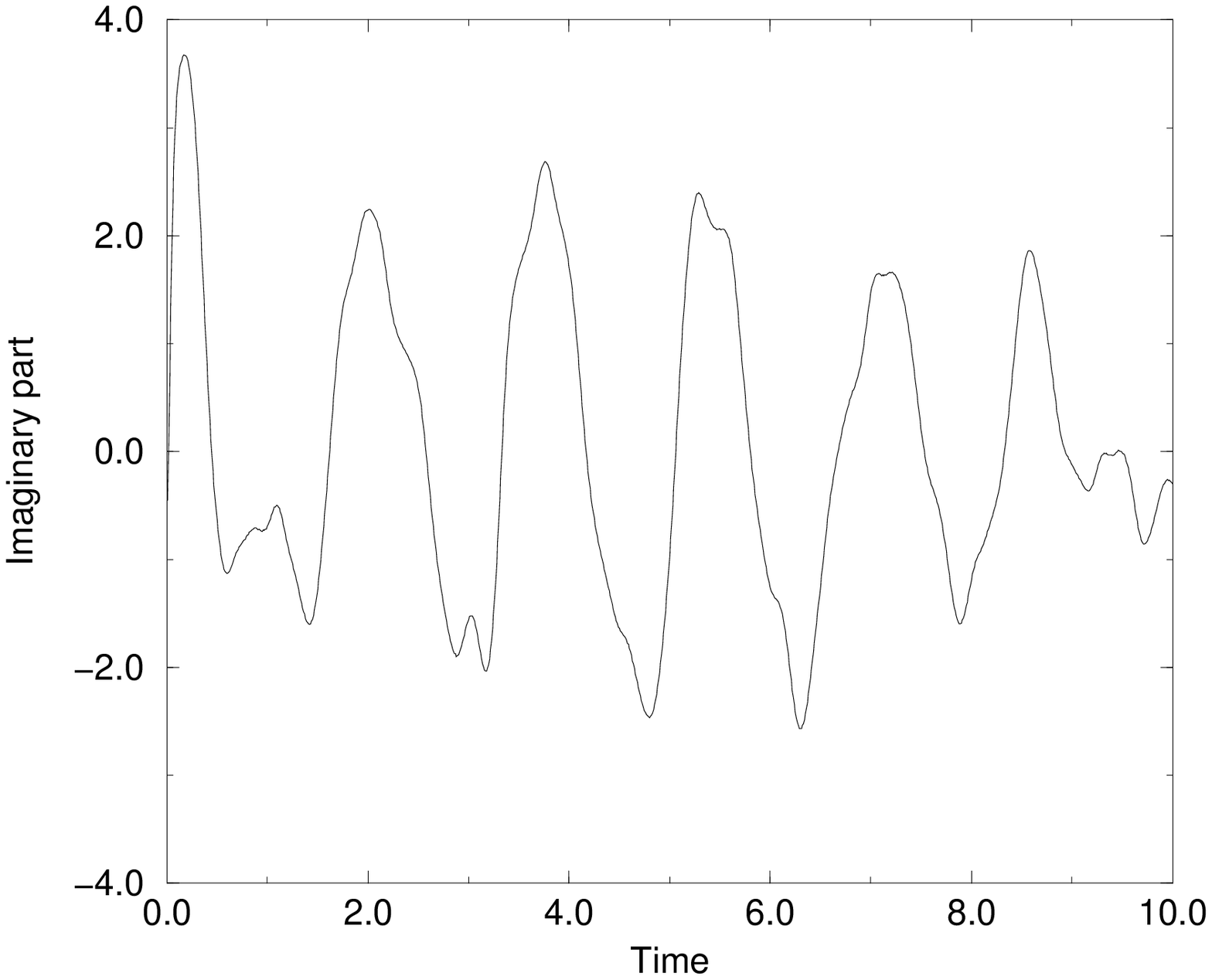,height=1.7in}}
\caption{Projection coefficient of the memory term for the
equation for the resolved mode 2 on the resolved mode 4 
for the second set of resolved variables. 
a) Real part, b) Imaginary part.}
\label{fig:res12}
\end{figure}

Fig.(\ref{fig:res12}) 
shows the typical behavior of the projection coefficients. As
is evident, the projection coefficients $(LQe^{sL}QLu_j,u_k)$ of
the short-memory approximation do not decay fast. There are still some 
coefficients that decay fast, but the majority of the projection coefficients 
do not and those determine the accuracy of the conditional expectation estimate. 
We only know (recall (\ref{sm41})) that the coefficients 
$(LQe^{sL}QLu_j,u_k)$ approximate well the quantities 
$(LQe^{sQL}QLu_j,u_k)$ for short times. Even though the coefficients 
$(LQe^{sL}QLu_j,u_k)$ of the short-memory approximation do 
not decay fast, the coefficients $(LQe^{sQL}QLu_j,u_k)$ can very well 
start decaying fast after a short time. However, we don't know that before 
implementing the short-time approximation. By the same reasoning as
in the case of resolution of all unstable modes, we have to truncate the interval 
of integration for the memory term, from $[0,t]$ to $[t-t_0,t]$ (for $t > t_0$). 
Again, $t_0$ is a short time over which the error of using 
$(LQe^{sL}QLu_j,u_k)$ instead of $(LQe^{sQL}QLu_j,u_k)$ is not large.  
Of course, we do not know the value of $t_0.$ In the case of resolution
of all unstable modes, we guessed that the interval $t_0$ should be
close to the time needed for the projection coefficients to decay significantly. 
But, in the present case, where the projection coefficients do not decay
fast, we cannot make any guess. We conducted numerical experiments with 
different truncated intervals of integration for the memory term and the results
do not change much.We present results for the case where the interval of
integration of the memory term is truncated from $[0,t]$ to $[t-1,t]$ (for
$t>1$).

The system of random integrodifferential equations is solved
repeatedly for different realizations of the noise. The initial
condition for the resolved modes $-3,-2,2,3$ is 
$u_2=u_{-2}=u_3=u_{-3}=1$, while the rest
of the resolved modes are set initially equal to zero.

\begin{figure}
\centering
\subfigure[]{\epsfig{file=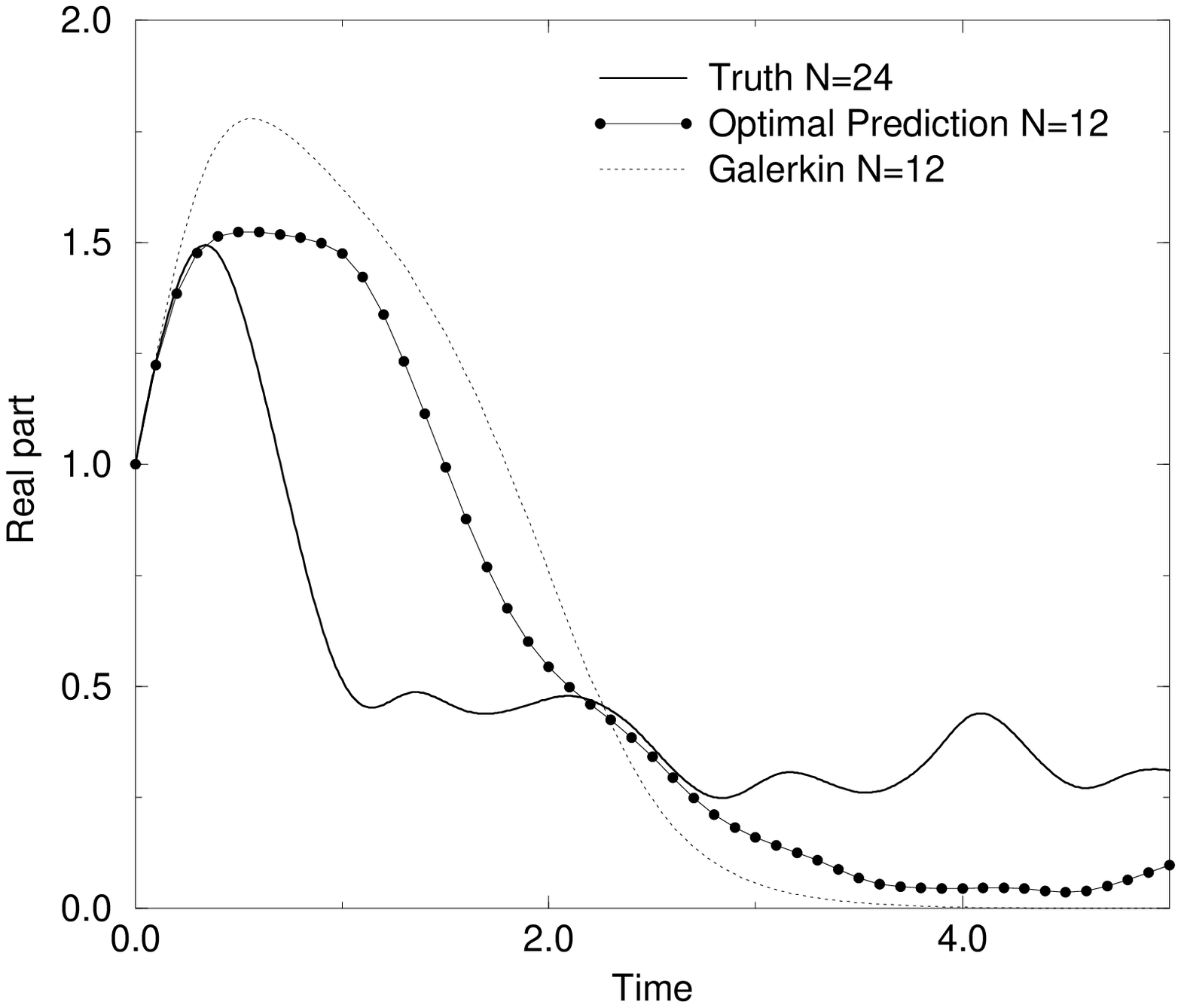,height=1.7in}}
\quad
\subfigure[]{\epsfig{file=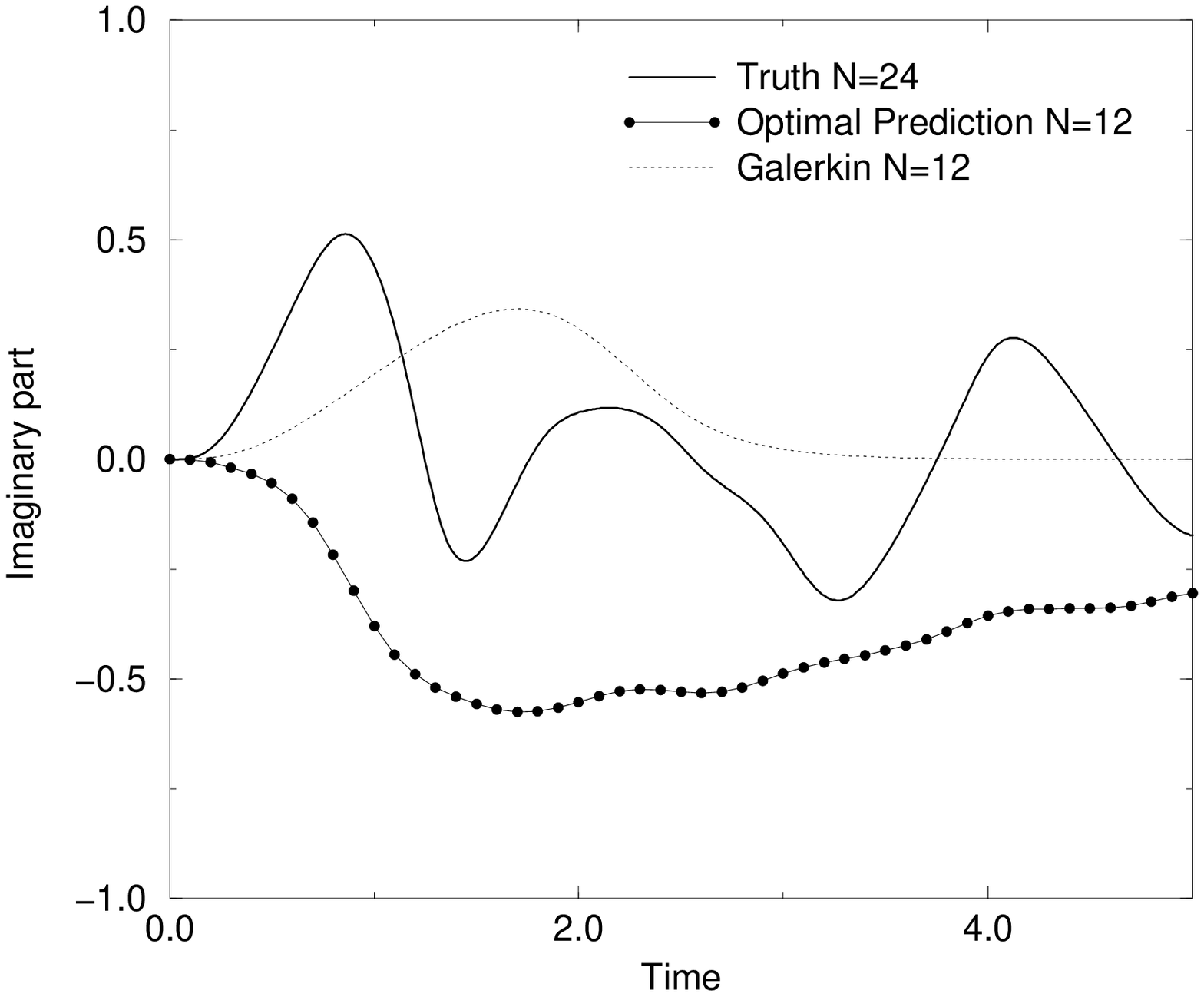,height=1.7in}}
\caption{Conditional expectation evolution for the
 resolved mode 2 for the second set of resolved variables. 
 a) Real part, b) Imaginary part.}
\label{fig:res13}
\end{figure}

\begin{figure}
\centering
\subfigure[]{\epsfig{file=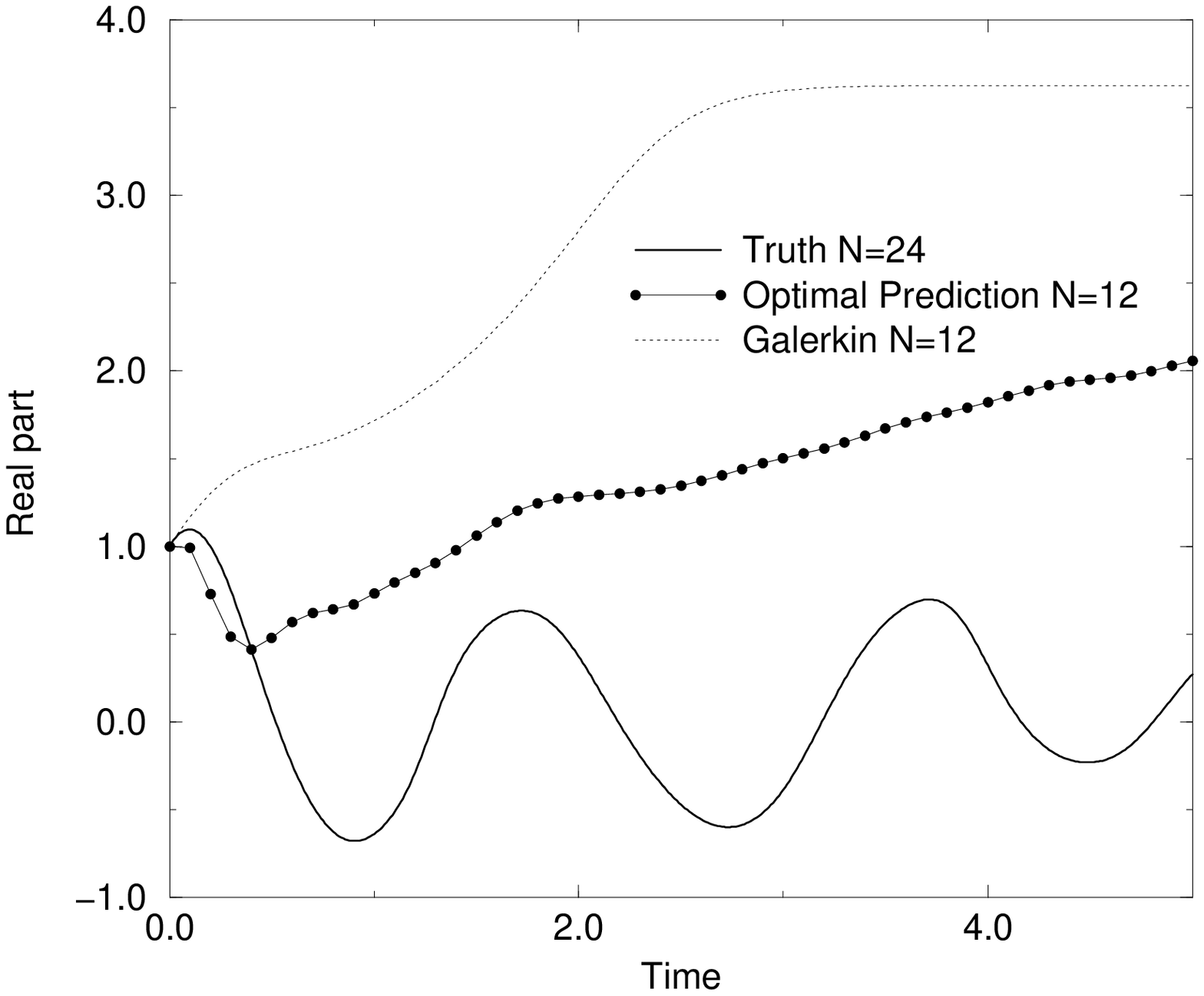,height=1.7in}}
\quad
\subfigure[]{\epsfig{file=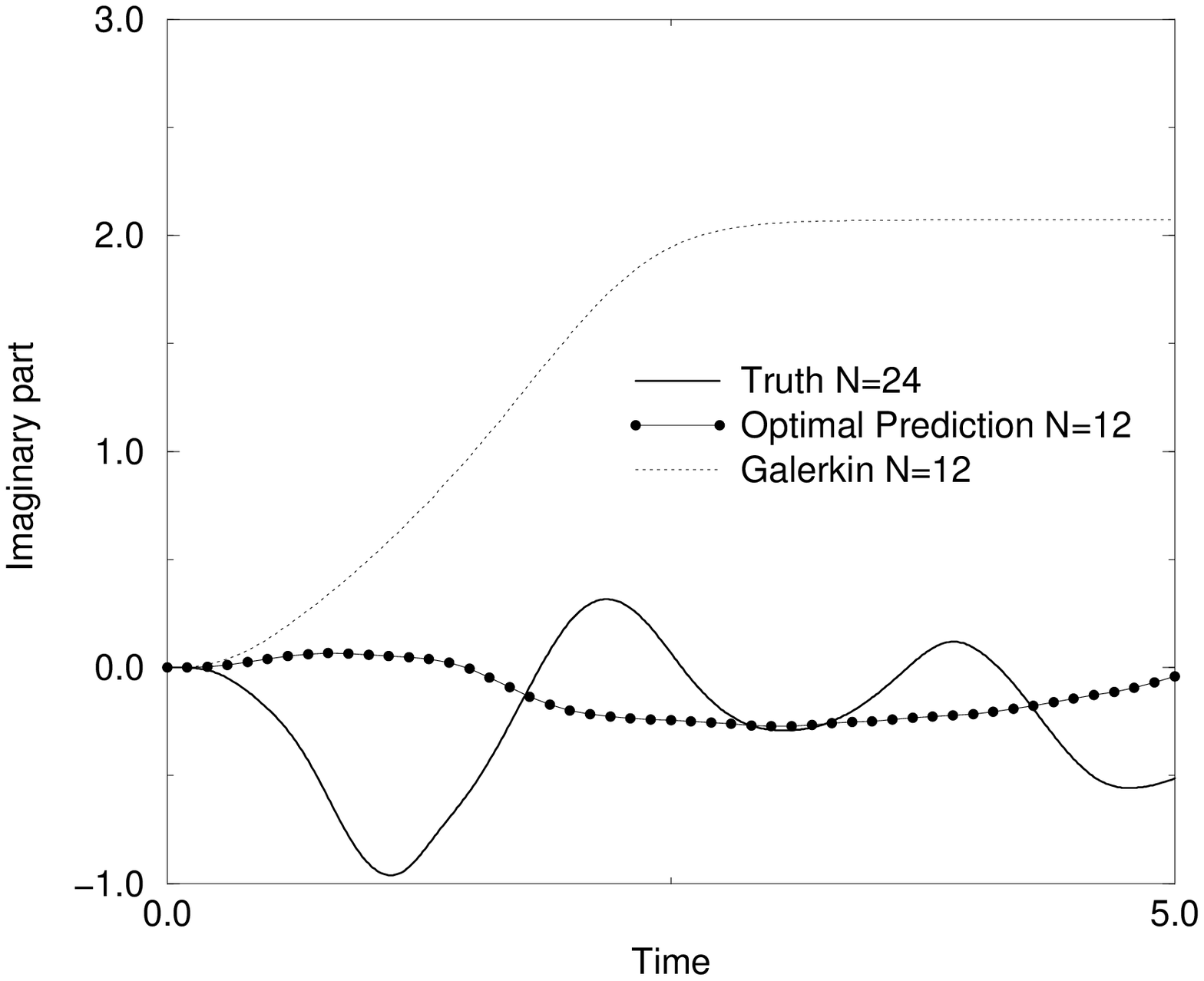,height=1.7in}}
\caption{Conditional expectation evolution for the
 resolved mode 3 for the second set of resolved variables. 
 a) Real part, b) Imaginary part.}
\label{fig:res14}
\end{figure}

\begin{figure}
\centering
\subfigure[]{\epsfig{file=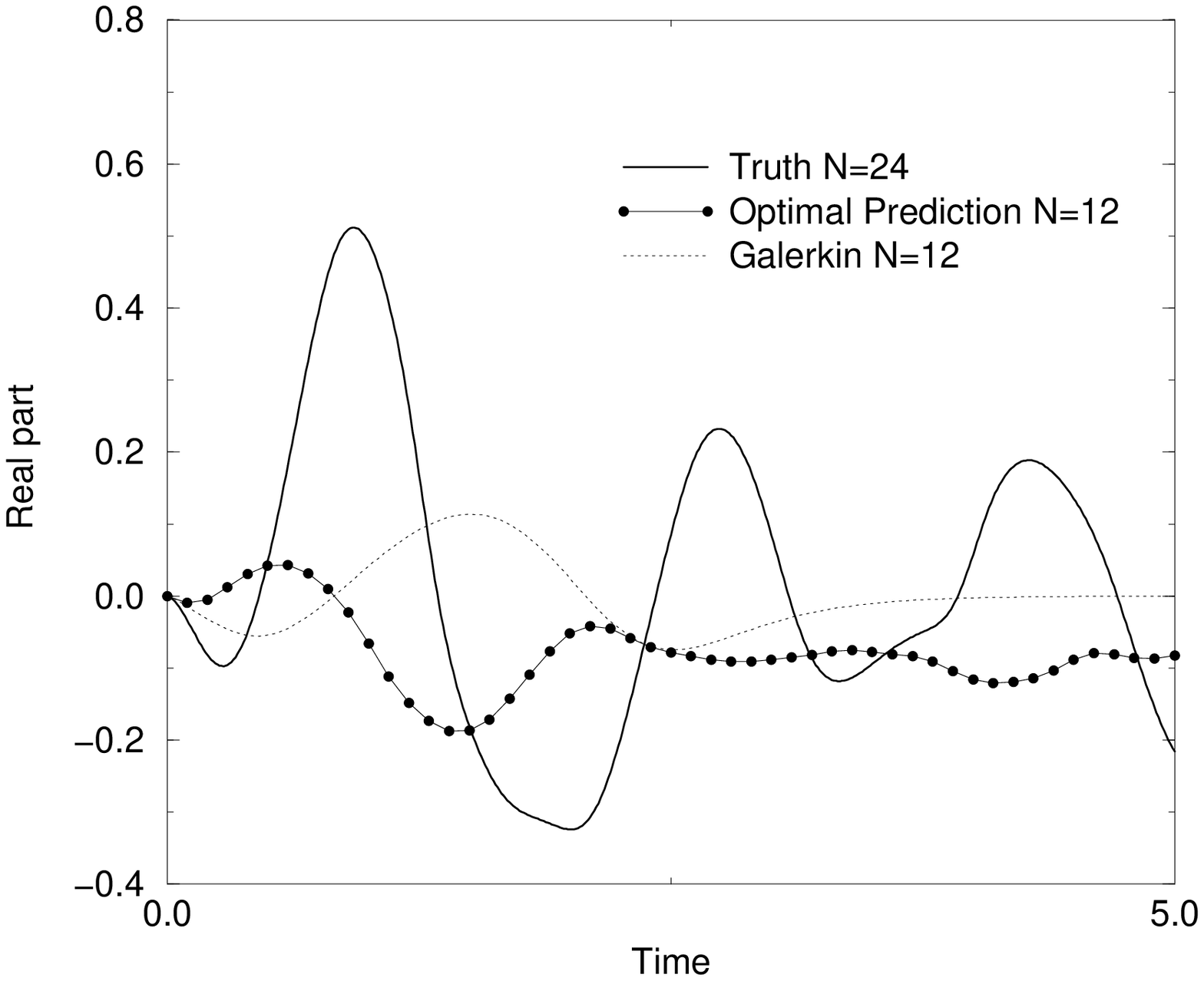,height=1.7in}}
\quad
\subfigure[]{\epsfig{file=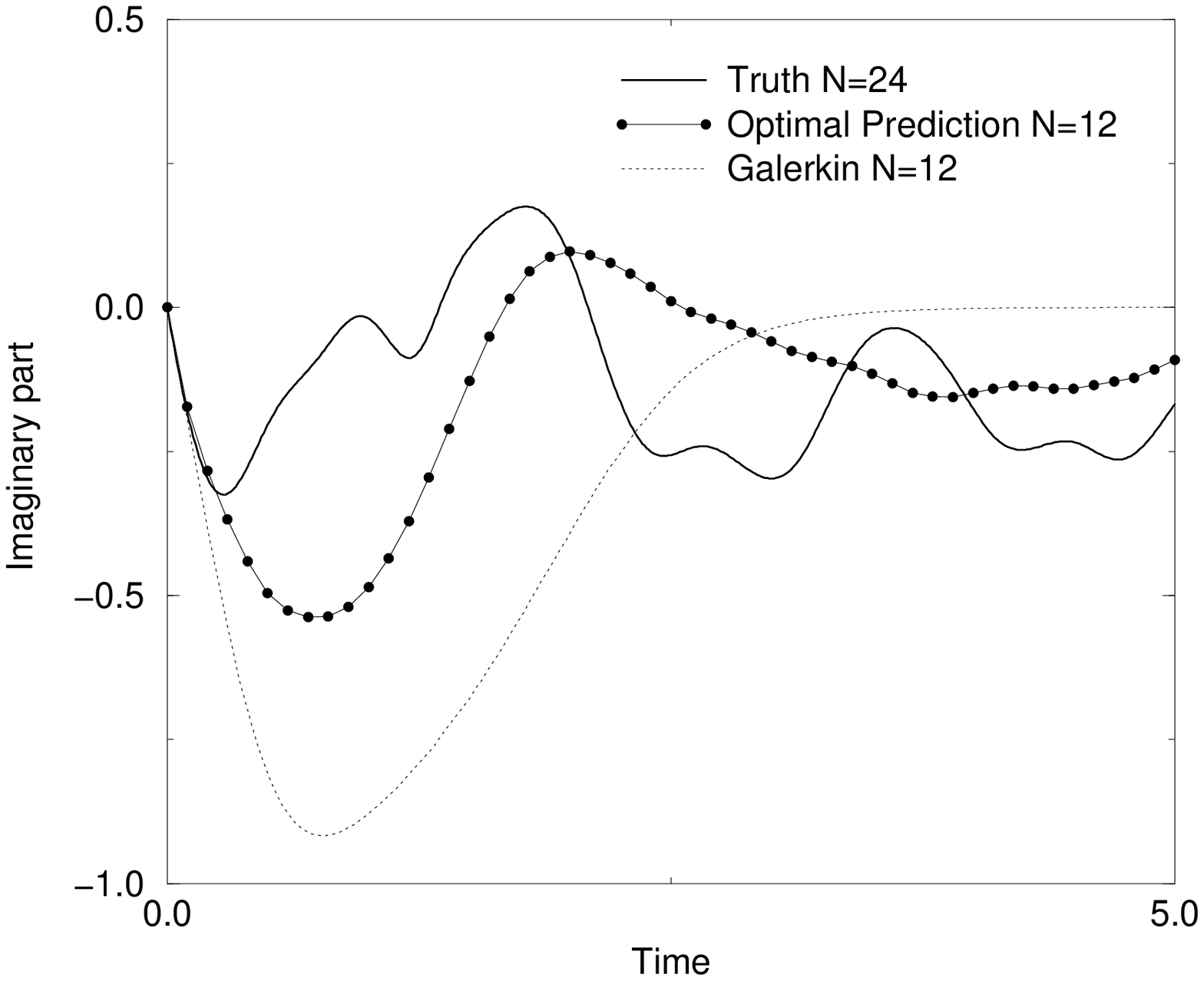,height=1.7in}}
\caption{Conditional expectation evolution for the
 resolved mode 4 for the second set of resolved variables. 
 a) Real part, b) Imaginary part.}
\label{fig:res15}
\end{figure}

\begin{figure}
\centering
\subfigure[]{\epsfig{file=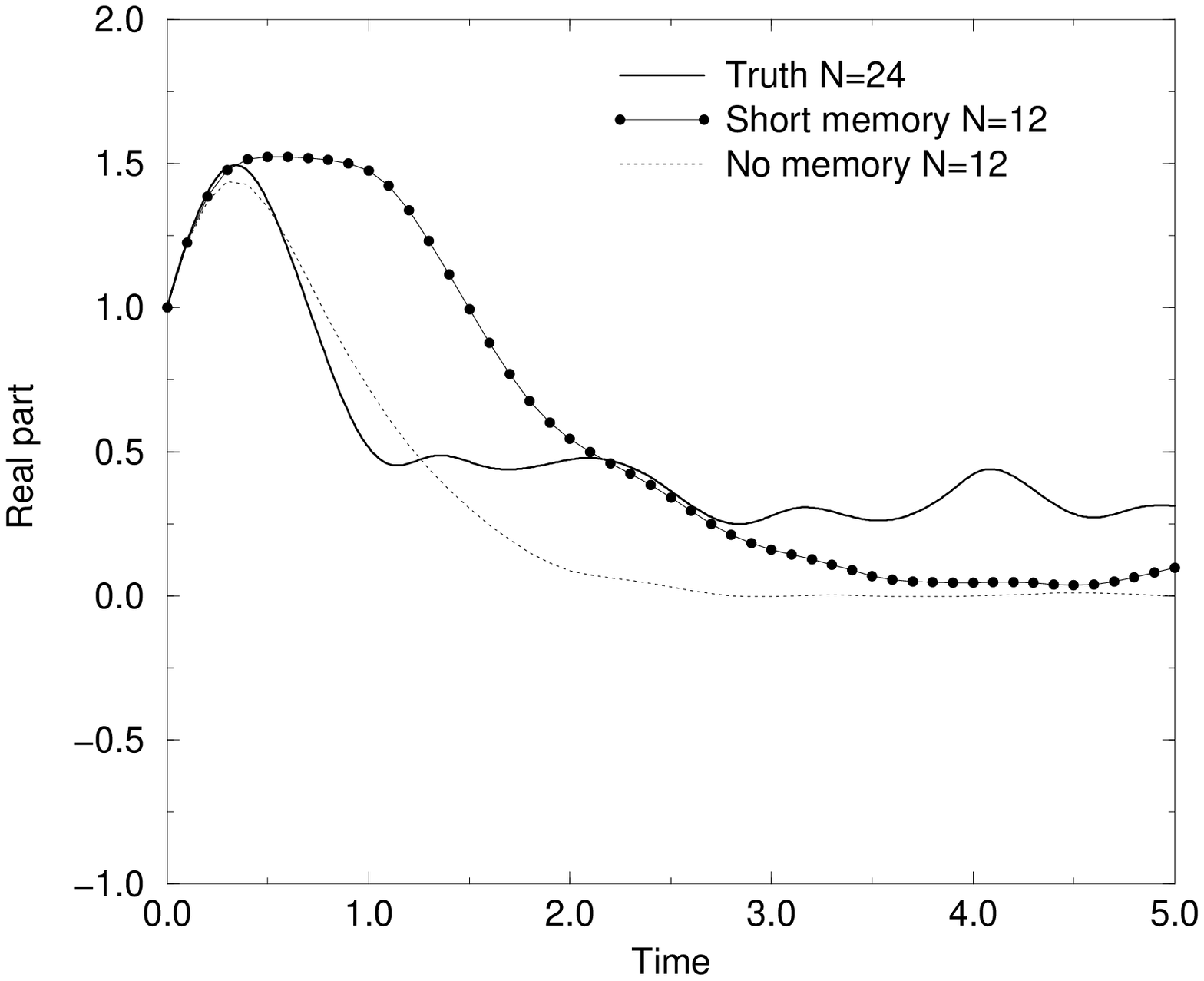,height=1.7in}}
\quad
\subfigure[]{\epsfig{file=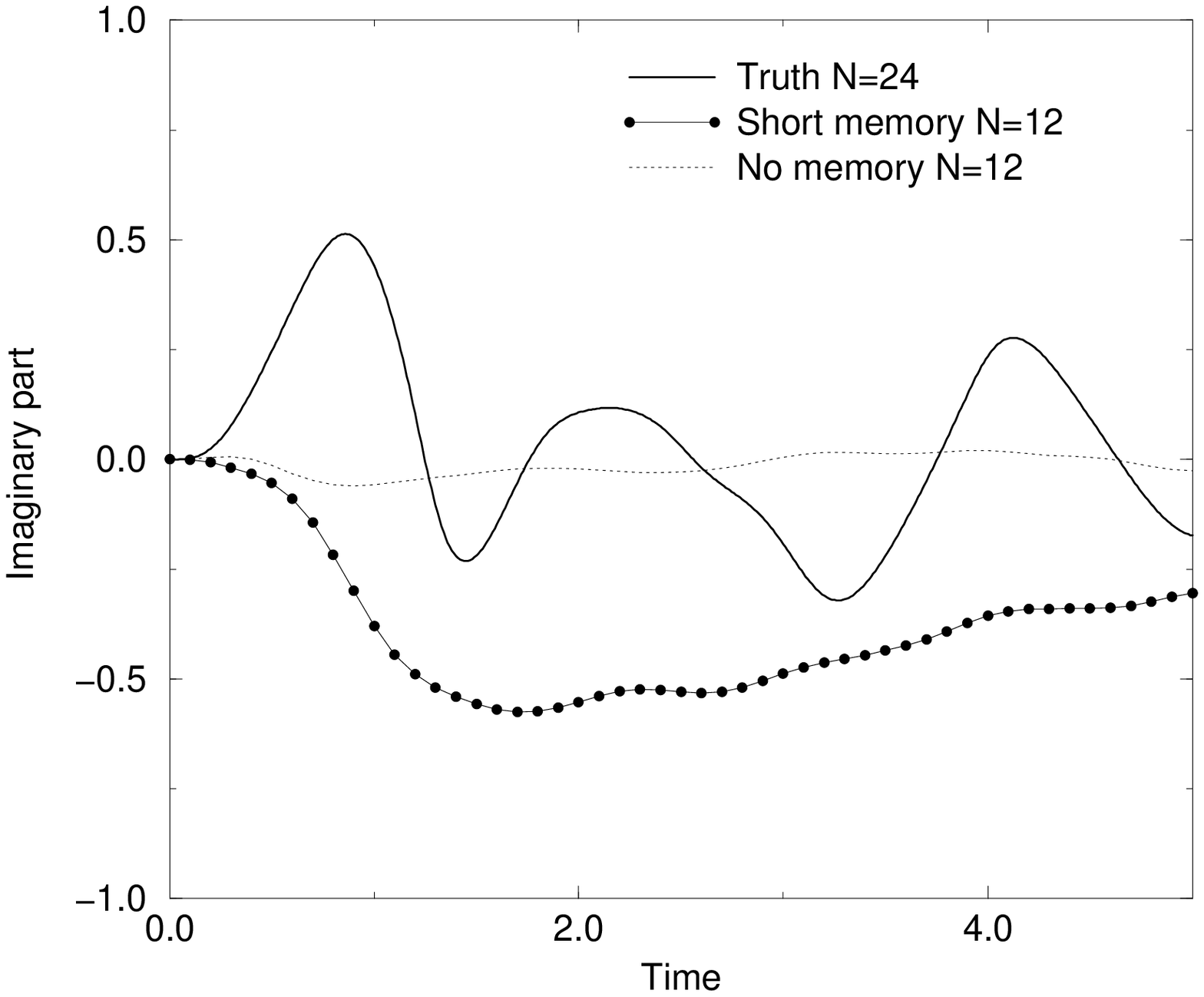,height=1.7in}}
\caption{Comparison of the conditional expectation evolution estimates for the
resolved mode 2, as produced by the short-memory approximation and
the short-memory approximation with the memory term set to zero.
a) Real part, b) Imaginary part.}
\label{fig:res13p}
\end{figure}

Figs.(\ref{fig:res13}), (\ref{fig:res14}), 
(\ref{fig:res15}) show the real and imaginary parts of the short-
memory estimates of the conditional expectations for the resolved modes 2,3,4
as compared to the truth and the Galerkin approximation. The
true conditional expectations were computed by averaging over
10000 samples. The short-memory estimates were computed
by averaging over 10000 realizations of the noise. As can be seen, 
there is good agreement between
the short-memory approximation estimates for the conditional
expectations of the resolved modes and the true conditional
expectations of these modes only for very short times. For longer
times the errors become large.

We show next that 
the results for the short-memory approximation 
depend mostly on the accuracy with which we approximate the memory 
term, and not so much on the accuracy of the approximation of the noise term. 
It is instructive to integrate the optimal prediction equations in the 
short-memory approximation with the memory term set to zero. 
Fig.(\ref{fig:res13p}) shows the conditional expectation estimate for the 
resolved mode 2 as predicted by the short-memory approximation when the 
memory term is set to zero. We compare it to the truth and the short-memory 
approximation when the memory term is not set to zero. As is
evident, the estimate is good for short times, while for longer times, the 
estimate converges to zero. Any deviations of the conditional expectations 
from zero, i.e. any improvement of the estimate, must be produced by the 
memory term. The deviations from zero, for long times, of the true 
conditional expectations, translate into memory effects for the set of 
resolved modes, and if these memory effects are not well represented 
in the model for the resolved modes, the prediction of the deviations is 
not accurate.

Since the short-memory approximation estimates are not good for long times, 
we conclude that the quantities $(LQe^{sQL}QLu_j,u_k)$ do not start decaying 
fast after a short time. If they did, the short-memory approximation would give 
good results for longer times. The inclusion of an unstable mode in the 
unresolved modes, results in the appearance of long-time memory effects, but 
the precise mechanism of formation of these long-time memory effects has
yet to be determined.

For the second set of variables, no attempt was made to compute the 
conditional expectation estimates using the finite-rank projection or the 
delta-function approximation. For the case of the finite-rank projection,
there is no chance of improvement, since the linear projection coefficients
that are part of the finite-rank projection, are already slowly-decaying.
For the case of the delta-function approximation, there is no incentive
to try it, since the projection coefficients do not decay fast.

%%%%%%%%%End of Section{Numerical simulations}%%%%%%

\section{Conclusions} 

We have applied the optimal prediction formalism to the solutions of the
Kuramoto-Sivashinsky equation in a Fourier-Galerkin truncation when 
some initial data are missing. The conditional expectations of the resolved 
modes conditioned on their initial values were estimated through simulation 
of the optimal prediction equations and compared to the true conditional 
expectations computed from simulations of the full system. The results of the 
comparison depend on which Fourier modes are contained in the set of resolved 
variables and on the type of projection used in the memory term. For the case 
when the resolved variables include all the unstable modes, we used two 
different kinds of projection for the memory term, namely the linear projection 
and the finite-rank one.

For the linear projection, the agreement between the optimal prediction 
estimates of the conditional expectations and the true conditional expectations is 
good for relatively long times. Also, the estimates show a considerable 
improvement over the Galerkin approximation, where we resolve a reduced set 
of variables and set all the unresolved variables equal to zero. For the case of the 
linear projection we, also, tried the more drastic delta-function approximation
where the correlations appearing in the memory term integrand are replaced by 
a delta-function multiplied by the integral. In this case, the resulting optimal
prediction equations are differential (not integrodifferential) equations. The 
results are not as good as in the more sophisticated short-memory 
approximation. However, the accuracy of the conditional expectation estimates
based on the delta-function approximation is impressive, if we consider
the much greater numerical efficiency of the resulting system of equations
for the resolved modes.

For the finite-rank 
projection, the agreement is good only for short times. The growth of the error 
for larger times is due to the fact that the short-memory approximation 
assumption is violated by the appearance of long-time memory effects. These
manifest themselves as slowly-decaying projection coefficients in the
memory term integral.

For the case where one unstable mode is left unresolved, we used only
the linear projection. The error of the estimates of the conditional expectations 
becomes large after a short time. The growth of the error is  due to the fact that 
the short-memory approximation assumption is violated by the appearance of 
long-time memory effects. The short-memory approximation estimates of the 
conditional expectations are still an improvement over the Galerkin 
approximation estimates. For the case where one unstable mode is left 
unresolved, we did not attempt to use the finite-rank projection, since the
linear projection that is part of the finite-rank projection, already violates 
the short-memory approximation.

In the case of 
the finite-rank projection for the first set of resolved modes, the long-time 
memory effects are a result of projecting on polynomials of the resolved modes
of order higher than 1. For the second set of resolved modes, 
where we used the linear projection for the memory term, the long-time 
memory effects are a result of leaving an unstable mode as unresolved. 
Although the long-time memory effects for these two cases have different 
causes, the result is the same, namely the violation of the short-memory 
approximation assumption. This leads to a large increase of the error in the 
conditional expectation estimates after a short time. The inadequacy of 
the short-memory approximation to produce accurate estimates of the
conditional expectations for long times, suggests that for cases with
slowly-decaying projection coefficients, the calculation of the quantities 
$(LQe^{sQL}QLu_j,u_k)$ cannot be avoided (see \cite{chorin1} for methods
of computing these quantities).

The determination of the reasons for the appearance of long-time memory
effects is important since it can provide insight on the way that the partial 
information known initially is lost when we perform underresolved 
computations. It can also help 
in determining which variables (or combinations of variables) of the full system 
should be resolved, if we hope to obtain a reduced model of the system that is 
accurate for long times \cite{chorin5,kreiss}.

%%%%%%%%%End of Section{Conclusions}%%%%%%%%%%%

\section*{Acknowledgments}
The work presented here is part of my doctoral thesis \cite{stinis}, and was conducted 
under the supervision of Professor Alexandre Chorin of UC Berkeley while I was a guest at the 
Lawrence Berkeley National Laboratory. I would like to thank Professor Chia-Kun Chu of 
Columbia University in New York for his generous support. Also, I would like to thank 
Dr. Kevin Lin and Dr. Mayya Tokman of UC Berkeley for very helpful discussions and comments.
 
\bibliographystyle{siam}

\bibliography{opks}

\end{document}